\documentclass[reqno, makeidx]{amsart}
        
\input supp-pdf.tex
\usepackage{graphicx}
\usepackage{latexsym}
\usepackage{amstext}
\usepackage {amsmath}
\usepackage {amsfonts}
\usepackage {amssymb}
\usepackage {amsthm}
\usepackage {bbm}
\usepackage{enumerate}
\usepackage{array}
\usepackage{bbm}
\usepackage[bookmarks=true]{hyperref}
\usepackage{stmaryrd}
\usepackage{mathrsfs}

\DeclareMathAlphabet{\mathpzc}{OT1}{pzc}{m}{it}

\DeclareMathOperator{\dist}{dist}

\newcommand{\res}      {\mathop{\hbox{\vrule height 7pt width .5pt depth
                        0pt\vrule height .5pt width 6pt depth
0pt}}\nolimits}

\newcommand{\ep}{\varepsilon}
\newcommand{\eps}{\varepsilon}

\newcommand{\grad}{\nabla}

\newcommand{\Ha}{\ensuremath{\mathcal{H}}}

\newcommand{\Hv}{\mathbf{H}}

\newcommand{\LL}{\ensuremath{\mathcal{L}}}

\def\muepsueps{\mu^\eps_{u_\eps}}
\newcommand{\N}{\ensuremath{\mathbb{N}}}
\def\niu{\nu_u}
\def\niueps{\nu_{u_\eps}}

\def\proju{P^u}
\def\projueps{P^{u_\eps}}

\def\R{\mathbb R}
\def\rectifiableset{M}
\newcommand{\Rn}{\ensuremath{{\mathbb{R}^n}}}

\newcommand{\sffu}{\mathbf{B}_{u}}
\newcommand{\sff}{\mathbf B}

\newcommand{\surftens}{c_0}

\def\tildemuepsueps{\widetilde\mu^\eps_{u_\eps}} 
 
\newcommand{\var}{\mathbf v}
 
\def\Vepsueps{V^\eps_{u_\eps}}

\newcommand{\mut}{\widetilde{\mu}}
\newcommand{\Vepstueps}{V^{\eps}_{\widetilde u_{\eps}}}
\newcommand{\Vepstuepsk}{V^{\eps_k}_{\widetilde u_{\eps_k}}}
\newcommand{\tuepsk}{\widetilde u_{\eps_k}}
\newcommand{\nituepsk}{\nu_{\tuepsk}}
\newcommand{\tueps}{\widetilde u_\eps}
\newcommand{\teps}{\widetilde\eps}
\newcommand{\tu}{\widetilde{u}}
\newcommand{\ind}{\mathbbm{1}}

\newcommand{\Galim}{\Gamma(L^1(\Omega))-\lim_{\eps\to 0}}
\newcommand{\Etildeps}{\widetilde{\mathcal E}_\eps}
\newcommand{\Eeps}{\mathcal E_\eps}

\theoremstyle{plain}
\numberwithin{equation}{section}
\newtheorem{lemma}{Lemma}[section]
\newtheorem{theorem}[lemma]{Theorem}
\newtheorem{proposition}[lemma]{Proposition}
\newtheorem{definition}[lemma]{Definition}

\newtheorem{corollary}[lemma]{Corollary}
\theoremstyle{definition}
\newtheorem{remark}[lemma]{Remark}
\begin{document}
\title{Gamma-convergence results for phase-field approximations of the 2D-Euler Elastica Functional}
\author{Luca Mugnai}
\address{Luca Mugnai, Max Planck Institute for Mathematics in the
  Sciences, Inselstr. 22, D-04103 Leipzig, Germany}
\email{mugnai@mis.mpg.de, }

\subjclass[2000]{Primary, 49J45; Secondary, 34K26, 49Q15, 49Q20}

\keywords{$\Gamma$-convergence, Relaxation, Singular Perturbation, Geometric Measure Theory}

\date{\today}
\begin{abstract}
We establish some new results about the $\Gamma$-limit, with respect to the $L^1$-topology, of two different (but related) phase-field approximations of the so-called Euler's Elastica Bending Energy for curves in the plane.
\end{abstract}

\maketitle

\pagestyle{myheadings}
\thispagestyle{plain}
\markboth{}{}
\section{Introduction}

In this paper we present some new results about the sharp interface limit of two sequences of phase-field functionals involving the  so-called Cahn-Hilliard energy functional and its $L^2$-gradient. The study of this kind of problems is motivated by applications in different fields ranging from image processing (\textit{e.g.}, \cite{EseShen, LoMar, Suka, BraMar}), to the diffuse interface approximation of elastic bending energies (\textit{e.g.}, \cite{DuWill, DuCaz, DuUno, Lowe, BeMu:09, CampHern2}), to the study of singular limits of partial differential equations and systems (\textit{e.g.}, \cite{Sato, MuRo2, LiuSatoTone, Serfa}), up to the study of rare events  for stochastic perturbations of  the so-called Allen-Cahn equation (\textit{e.g.}, \cite{KORV, MuRo1}).  
 
Let us now introduce the two  sequences of energies we wish to study. Given  $\Omega\subset\R^d$ open, bounded and with smooth boundary, we define the so-called Cahn-Hilliard energy by 
\begin{gather}\label{eq:def-Van-der-Waals}
\mathcal P_\eps(u):=\begin{cases} \int_{\Omega}\frac{\eps}{2}\vert\nabla u\vert^2+\frac{W(u)}{\eps}\,dx & \text{if }u\in W^{1,2}(\Omega),
\\
+\infty & \text{otherwise on }L^1(\Omega)
\end{cases}
\end{gather}
where  $\eps>0$ is a parameter representing the typical ``diffuse interface width'', and $W\in C^3(\R,\R^+\cup\{0\})$ is a double-well potential with two equal minima (throughout the paper we make the choice $W(s):=(1-s^2)^2/4$, though most of the results we obtain hold true for a wider class of potentials).
The sequences of functionals  $\{\Etildeps\}_\eps,\,\{\Eeps\}_\eps$ we consider in this paper are respectively defined  by 
\begin{gather}
\widetilde{\mathcal E}_{\eps}:= \mathcal P_\eps+\mathcal W_\eps
:\,L^1(\Omega)\to\,[0,+\infty],
\\
\text{ where  }
\mathcal W_\eps(u):=\begin{cases}\frac{1}{\eps}\int_{\Omega}\Big(\eps\Delta u-\frac{W^\prime(u)}{\eps}\Big)^2\,dx & \text{if }u\in W^{2,2}(\Omega)
\\
+\infty &\text{elsewhere on }L^1(\Omega),
\end{cases},
\label{eq:def-func-DG}
\end{gather}
and 
\begin{gather}
\mathcal E_{\eps}:= \mathcal P_\eps+\mathcal B_\eps
:\,L^1(\Omega)\to\,[0,+\infty],
\\
\text{where }
\mathcal B_\eps(u):=\begin{cases}\frac{1}{\eps}\int_{\Omega}\Big\vert\eps\nabla^2 u-\frac{W^\prime(u)}{\eps}\nu_{u}\otimes\nu_{u}\Big\vert^2\,dx & \text{if }u\in W^{2,2}(\Omega)
\\
+\infty &\text{elsewhere on }L^1(\Omega)
\end{cases},
\end{gather}
and $\nu_{u}$ is a unit vector-field  such that
$$ 
\nu_{u}=\frac{\nabla u}{\vert\nabla u\vert}\,\text{on }
\{\nabla u\neq 0\}\text{ and } \niu\equiv const.\,\text{on }\{\nabla u= 0\}.
$$
We remark that $\mathcal W_\eps(u)$ represents the (rescaled) norm of the $L^2$-gradient of $\mathcal P_\eps$ at $u$, and that  $\mathcal W_\eps$ and $\mathcal B_\eps$ are linked by the relation
$$
\mathrm{tr}\Big[\eps\nabla^2 u
-\frac{W^\prime(u)}{\eps}\nu_u\otimes\nu_u\Big]=\eps\Delta u-\frac{W^\prime(u)}{\eps}.
$$
Hence, denoted by $\{\lambda_1,\dots,\lambda_d\}$ the eigenvalues of the symmetric $d\times d$-matrix $\eps\nabla^2u-(W^\prime(u)/\eps)\nu_u\otimes\nu_u$, we have
\begin{align}
d\Big(\sum_{i=1}^d\lambda_i^2\Big)=d\left\vert\eps\nabla^2 u
-\frac{W^\prime(u)}{\eps}\nu_u\otimes\nu_u\right\vert^2
\geq 
\left(\eps\Delta u_\eps-\frac{W^\prime(u_\eps)}{\eps}\right)^2
=\Big(\sum_{i=1}^d\lambda_i\Big)^2.
\label{smanettando}
\end{align}

Next, we briefly summarize the known results about the sharp interface limit of $\{\widetilde{\mathcal E}_\eps\}_\eps$ and $\{\mathcal E_\eps\}_\eps$. The starting point for the analysis of the asymptotic behavior, as $\eps\to 0$, of the sequences $\{\widetilde{\mathcal E}_\eps\}_\eps,\, \{\mathcal E_\eps\}_\eps$ is  a well-known result, due to Modica and Mortola, establishing  the $\Gamma$-convergence of $\mathcal P_\eps$ to the area functional. More precisely in  \cite{MM} it has been proved that the $\Gamma(L^1(\Omega))$-limit  of the sequence $\{\mathcal P_\eps\}_\eps$  is given by 
\begin{gather*}
\Gamma(L^1(\Omega))-\lim_{\eps\to 0}\mathcal P_\eps(u)=\mathcal P(u):=\begin{cases}\frac{\surftens}{2}\int_\Omega\,d\vert\nabla u\vert
& \text{if }u\in BV(\Omega,\{-1,1\}),
\\
+\infty & \text{elsewhere in }L^1(\Omega)
\end{cases}
\end{gather*}
where   $\surftens:=\int_{-1}^1\sqrt{2W(s)}\,ds$  (see Section \ref{bomboloni} and Section \ref{sec:none} for further details). We remark that for every $u\in BV(\Omega,\{-1,1\})$  we can write $u=2\chi_E-1=:\mathbbm{1}_E$, where $\chi_E$ denotes the characteristic function of the finite perimeter set $E:=\{u\geq 1\}$. Hence $\mathcal P(u)=\surftens\Ha^{d-1}(\partial^*E)$ where $\Ha^{d-1}$ denotes the $(d-1)$-dimensional Hausdorff measure in $\R^d$ and $\partial^* E$ denotes the reduced boundary of $E$ (see \cite{Sim}).  

The main result concerning the  $\Gamma$-convergence of $\{\widetilde{\mathcal E}_\eps\}_\eps$ has been established, for $d=2$ and $d=3$, by R\"oger and Sch\"atzle in \cite{RoSch} and independently, but only in the case $d=2$, by Tonegawa and Yuko in \cite{ToneYuko}, partially answering to a conjecture of De Giorgi (see \cite{DG}).  In particular in \cite{RoSch} the authors  proved  that for $d=2$ or $3$ and $u=\ind_E\in BV(\Omega,\{-1,1\})$ such that  $E\subset \Omega$ is open and $\Omega\cap\partial E\in C^2$, we have
\begin{equation}\label{eq:RS-intro}
\Gamma(L^1(\Omega))-\lim_{\eps\to 0}\widetilde{\mathcal E}_\eps(u)=\surftens\int_{\Omega\cap\partial E}\big[1+\vert\Hv_{\partial E}(x)\vert^2]\,d\Ha^{d-1}(x),
\end{equation}
where $\Hv_{\partial E}(x)$ denotes the mean curvature vector of $\partial E$ in the point $x\in\partial E$.  When $d=2$ we call the functional on the right hand side of \eqref{eq:RS-intro} the Euler's Elastica Functional.

The sequence of functionals $\{\Eeps\}_\eps$ has been introduced in \cite{BeMu:09} in connection with the problem of finding a diffuse interface approximation of the Gaussian curvature. As a straightforward consequence of the results established in \cite{BeMu:09} it follows that, again for $d=2,3$ and  $u=\ind_E\in BV(\Omega,\{-1,1\})$  such that $E\subset\Omega$ is open and $\Omega\cap\partial E\in C^2$, we have
\begin{equation}\label{eq:BeMu-intro}
\Gamma(L^1(\Omega))-\lim_{\eps\to 0}\mathcal E_\eps(u)=\surftens\int_{\Omega\cap\partial E}[1+\vert \sff_{\partial E}(x)\vert^2]\,d\Ha^{d-1}(x),
\end{equation}
where this time $\sff_{\partial E}(x)$ denotes the second fundamental form of $\partial E$ in the point $x\in\partial E$. 

In the present paper we restrict to the case $d=2$, and investigate the behavior of $\{\widetilde{\mathcal E}_\eps\}_\eps$ and $\{\mathcal E_\eps\}_\eps$ along sequences $\{u_\eps\}_\eps\subset C^2(\Omega)$ such that
\begin{equation}\label{fregnacce}
L^1(\Omega)-\lim_{\eps\to 0}u_\eps= \ind_E\in BV(\Omega,\{-1,1\}),
\end{equation}
 removing the regularity assumption on the limit set $E$. In other words we aim to prove a full $\Gamma$-convergence result, on the whole space $L^1(\Omega)$. 

\noindent We recall that if a sequence of functionals $\Gamma$-converges, and a certain equicoercivity property holds, then the minimizers of such sequence converge to the minimizers of the $\Gamma$-limit. Therefore, besides its possible mathematical interest, we expect that a description of the $\Gamma$-limit may be of some relevance at least for those applications, such as \cite{BraMar, MeMaPa,DuWill, DuCaz, 
 BeMu:09, Lowe}, where  the sequences $\{\widetilde{\mathcal E}_\eps\}_\eps$ and $\{\mathcal E_\eps\}_\eps$ are introduced to formulate, and solve numerically, a ``diffuse interface'' variational problem whose solutions are expected to converge, as $\eps\to 0$, to the solutions of a given sharp interface minimum problem.

In synthesis our result says that the sharp interface limits of  $\{\Etildeps\}_\eps$ and $\{\Eeps\}_\eps$ in general do not coincide out of ``smooth points'', although in two space dimensions, by  \eqref{eq:RS-intro} and \eqref{eq:BeMu-intro}
and 
 \begin{gather}\label{eq:2d}
\vert \sff_{\partial E}(x)\vert^2=\vert \Hv_{\partial E}(x)\vert^2,
\end{gather}
we have, for $u=\ind_E$ and $E\subset\Omega$  open such that  $\Omega\cap\partial E\in C^2$,
\begin{equation*}
\Galim\mathcal{E}_\eps(u)=
\Galim\widetilde{\mathcal{E}}_\eps(u).
\end{equation*}
 More precisely we prove that, in accordance with \eqref{smanettando}, we have
$$
\Galim\Etildeps(u)\leq
\Galim\Eeps(u)
$$
and that there are functions $u\in BV(\Omega,\{-1,1\})$ for which the above inequality is strict. In fact we show that, on one hand, a uniform bound on $\Eeps(u_\eps)$ implies that the  energy density measures $$
\mu_\eps:=\surftens^{-1}[\eps/2\vert\nabla u_\eps\vert^2+W(u_\eps)/\eps]\LL^d_{\res\Omega}
$$ 
(here $\LL^d$ denotes the Lebesgue measure on $\R^d$) concentrate on a set whose tangent cone in \textit{every} point is given by an unique tangent line. On the other hand we show that a uniform bound on $\Etildeps(u_\eps)$ allows the energy measures to concentrate on cross-shaped sets. This difference in regularity between the  support of the two limit measures is related to the existence of so called ``saddle shaped solutions'' to the semilinear elliptic equation $-\Delta U+W^\prime(U)=0$ on $\R^2$ (see \cite{Dang-Fife-Pel, CabTe:JEMS, DelPino}, and the proof of Theorem \ref{theo:sega} in this paper).

To give a better description of our results, let us briefly explain the role played by the regularity assumption on the limit set $E$ in the proofs of \eqref{eq:RS-intro} and \eqref{eq:BeMu-intro}, and discuss the obstructions to remove such an assumption. To this aim, 
for the readers convenience, we briefly recall the backbone of the proof of \eqref{eq:RS-intro} (and point out that the proof \eqref{eq:BeMu-intro} follows the same line of arguments).

To prove  \eqref{eq:RS-intro} one has to find a lower-bound for $\widetilde{\mathcal E}_\eps(u_\eps)$ proving the so-called $\Gamma-\liminf$ inequality; and to show that such lower-bound is in a way ``optimal'' via the so-called $\Gamma-\limsup$ inequality. (See Section \ref{bomboloni} for a precise definition of $\Gamma$-convergence).

\noindent Let us begin recalling how the $\Gamma-\liminf$ inequality has been proved.
Suppose that $\{u_\eps\}_\eps\subset C^2(\Omega)$ verifies \eqref{fregnacce} and 
\begin{gather}
\sup_{\eps>0}\widetilde{\mathcal E}_\eps(u_\eps)<+\infty.
\label{pippopazzo}
\end{gather} 
Thanks to the bound (uniform in $\eps$) on $\mathcal P_\eps(u_\eps)$ and the convergence of $u_\eps$ to $u$, applying the results of \cite{MM} it can be easily deduced that (up to subsequences) the  energy-density measures 
$
\mu_\eps$ defined above converge to a Radon measure $\mu$ in $\Omega$
such that
$
\Ha^{d-1}_{\res\partial^*E}<<\mu.
$
That is, roughly speaking, the support of $\mu$ contains the (reduced) boundary of $E$. In case only a bound on $\mathcal P_\eps(u_\eps)$ is available, there is no much hope to obtain more informations about the measure $\mu$, since this latter may be quite irregular (for example  it may contain parts that are absolutely continuous with respect to $\LL^d$). However when \eqref{pippopazzo} holds, the bound on $\mathcal W_\eps(u_\eps)$ implies that $\mu$ has some ``weak'' regularity properties. In fact 
the first crucial step in the derivation of a lower bound for $\widetilde{\mathcal E}_\eps(u_\eps)$ consists in proving that
\eqref{pippopazzo} guarantees that:
\begin{itemize}
\item the measure $\mu$ has the form  $\mu=\theta\,\Ha^{d-1}_{\res M}$,  where $M$ is a generalized hypersurface of $\R^d$, and $\theta:M\to\N$  is an integer valued  $\Ha^{d-1}_{\res M}$ -measurable function;
\item  
 a generalized mean curvature vector $\Hv_{\mu}\in L^2(\mu)$ is well defined $\mu$-a.e.; 
\item the following relation holds
\begin{equation}\label{eq:lb}
\liminf_{\eps\to 0} \widetilde{\mathcal E}_\eps(u_\eps)\geq \surftens\int[1+\vert\Hv_\mu\vert^2]\,d\mu.
\end{equation}
\end{itemize}
(Namely, $\mu$ is the weight measure of an integral rectifiable varifold with  $L^2$-bounded first variation, see Section \ref{luponeassassino} for some basic facts and terminology  about varifolds theory).

\noindent The next step in the proof of the $\Gamma-\liminf$-inequality consists in relating the generalized mean curvature vector $\Hv_\mu$ of $\mu$ with the (generalized) mean curvature vector of the phase-boundary $\partial^*E$. Since $\Ha^{d-1}_{\res\partial^*E}<<\mu$, by the results of \cite{Men} (see also \cite{Sch:09, Masnou}) it follows that  $\partial^*E$  can be covered with the union of a countable family of $(d-1)$-dimensional  $C^2$-manifolds embedded in $\R^d$, and with a set of $\Ha^{d-1}$-measure zero. Hence the mean curvature vector $\Hv_{\partial^*E}$ of $\partial^* E$ is well defined $\Ha^{d-1}_{\res\partial^*E}$-a.e. Furthermore by  \cite{Men} (see also \cite{Masnou, Sch:09} and Remark \ref{rem:local}) we have  $\Hv_\mu(x)=\Hv_{\partial^*E}(x)$ for $\Ha^{d-1}$-a.e. $x\in\partial^*E$ . Eventually, being $\Ha^{d-1}_{\res\partial^*E}<<\mu$ and $\theta(x)\geq 1$ for $\mu$-a.e. $x\in M$, by \eqref{eq:lb} it follows that
\begin{equation}\label{eq:ga-inf-intro}
\liminf_{\eps\to 0} \widetilde{\mathcal E}_\eps(u_\eps)\geq \surftens\int[1+\vert\Hv_\mu\vert^2]\,d\mu\geq \surftens\int_{\Omega\cap\partial^*E}[1+\vert\Hv_{\partial^* E}\vert^2]\,d\Ha^{d-1}.
\end{equation}
It remains to establish  if (or when) such a lower bound is ``optimal''. More precisely, it remains to understand for which $u=\ind_E\in BV(\Omega,\{-1,1\})$ is it  possible to find  a `` recovery sequence'' $\{u_\eps\}_\eps\subset C^2(\Omega)$, that is a sequence such that 
\begin{equation}\label{eq:ga-sup-intro}
\lim_{\eps\to 0}u_\eps=u\text{ in }L^1(\Omega)\text{ and } \limsup_{\eps\to 0}\widetilde{\mathcal E}_\eps(u_\eps)\leq \surftens\int_{\Omega\cap\partial^*E}[1+\vert\Hv_{\partial^* E}\vert^2]\,d\Ha^{d-1}.
\end{equation}
For those $u$ it follows that 
\begin{equation}\label{lamb}
\Gamma(L^1(\Omega))-\lim_{\eps\to 0}\widetilde{\mathcal E}_\eps(u_\eps)= \surftens\int_{\partial^* E\cap\Omega}[1+\vert\Hv_{\partial^* E}\vert^2]\,d\Ha^{d-1}
\end{equation}
(in fact \eqref{eq:ga-inf-intro} and \eqref{eq:ga-sup-intro} respectively represent the $\Gamma-\liminf$ and the $\Gamma-\limsup$-inequality). When $E\subset\Omega$ and $\Omega\cap\partial E\in C^2$, it is relatively easy to construct a sequence $\{u_\eps\}_\eps\subset C^2(\Omega)$ verifying \eqref{eq:ga-sup-intro} (see \cite{BePa:93}), and this concludes the proof of \eqref{eq:RS-intro}. 

Actually, by a simple diagonal argument, we can construct a sequence $\{u_\eps\}_\eps\subset C^2(\Omega)$ verifying \eqref{eq:ga-sup-intro} (and consequently obtain that \eqref{lamb} holds true) for all those functions $u=\ind_E\in BV(\Omega,\{-1,1\})$ for which there exists a sequence $\{E_h\}_{h}$ such that $E_h\subset\Omega$ is open and $\Omega\cap\partial E_h\in C^2$ for every $h\in\N$, and such that
$$
\lim_{h\to\infty}  \ind_{E_h}=u\text{ in }L^1(\Omega),\quad \lim_{h\to\infty}
\int_{\Omega\cap\partial E_h}[1+\vert\Hv_{\partial E_h}\vert^2]\,d\Ha^{d-1}=\int_{\Omega\cap\partial E}[1+\vert\Hv_{\partial E}\vert^2]\,d\Ha^{d-1}
$$
(\textit{e.g.}, if $E\subset\Omega$ is open and $\Omega\cap\partial E\in W^{2,2}$, see Remark \ref{rem:nadiacassini}). 

As we already said  \eqref{eq:BeMu-intro} is obtained following a similar line of arguments.

\begin{center}
 \begin{figure}[h]
      \includegraphics[width=.6\linewidth]{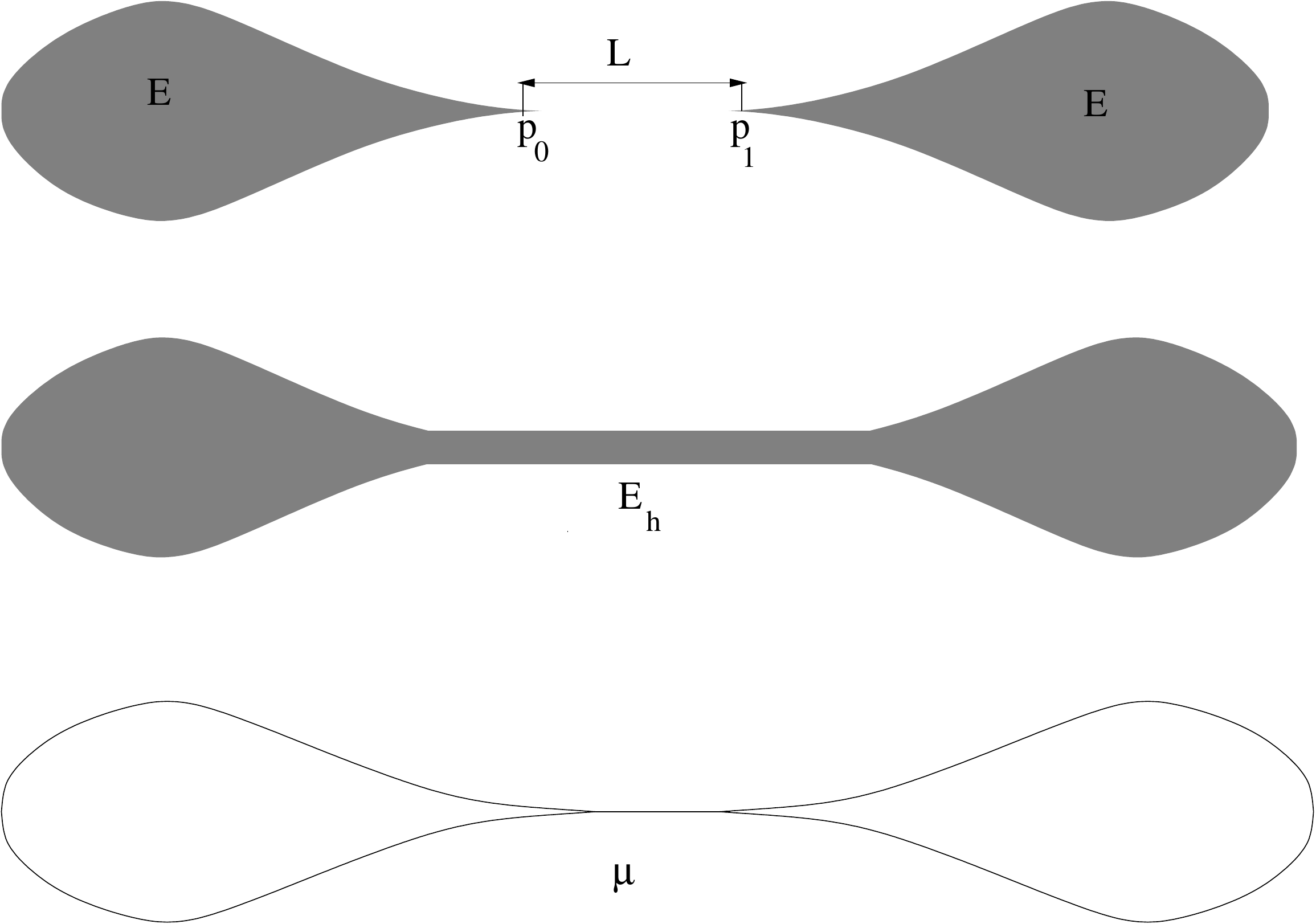} 
      \caption{\small The set $E\subset\subset\Omega$ has smooth boundary out of the two sharp cusps $p_0,\,p_1$, that are aligne and placed at a distance $L>0$. The first variation of the rectifiable varifold $\widehat V=\var(\partial^*E,1)$  associated to $\widehat\mu:=\Ha^{1}_{\res \partial^* E}$  is given by $\delta \widehat V=\Hv_{\partial^*E}\Ha^{1}_{\res \partial^* E}+2\sum_{j=0,1}\mathbf{e}_1(-1)^j\delta_{p_j}$, where $\mathbf e_1=(1,0)\in\R^2$ and $\delta_{p_j}$ is the Dirac-delta at $p_j$. Hence $\delta \widehat V\notin [L^2(\widehat\mu)]^*$, see Section \ref{luponeassassino}. For every $h\in\N$ the set $E_h$, such that $E_h\subset\subset\Omega$ and $\partial E_h\in C^2$,  is obtained replacing $p_0$ and $p_1$ with a flat tubular neighborhood of height $1/h$ that approximates the segment connecting $p_0$ and $p_1$ as $h\to\infty$.}
      \label{cas}
  \end{figure}
\end{center}
Yet we do not expect neither \eqref{lamb} nor its analogue for the $\Gamma(L^1(\Omega))-\lim_{\eps\to 0}\mathcal E_\eps(u)$ to be always true, as the following example suggests. Suppose that $E\subset\subset\Omega$ is as in Figure \ref{cas}. We then have $u=\ind_E\in BV(\Omega,\{-1,1\})$. Moreover if we consider the sequence $\{E_h\}_{h}$ of smooth sets represented in Figure \ref{cas}, we have
\begin{gather*}
L^1(\Omega)-\lim_{h\to\infty}\ind_{E_h}=u,
\\
\lim_{h\to\infty}\int_{\partial E_h}[1+\vert\Hv_{\partial E_h}\vert^2]\,d\Ha^1=\int[1+\vert\Hv_{\mu}\vert^2]\,d\mu=
\int_{\partial^* E}[1+\vert\Hv_{\partial^* E}\vert^2]\,d\Ha^1+2L.
\end{gather*}
Hence, for every $u_h:=\ind_{E_h}$ ($h\in\N$), by \cite{BePa:93, BeMu:09}, we can construct a recovery sequence $\{u_{h,\eps}\}_\eps\subset C^2(\Omega)$. Then, by a diagonal argument, we can select a sequence 
$\{u_\eps\}_\eps\subset C^2(\Omega)$ such that $\lim_{\eps\to 0} u_\eps=u$ in $L^1(\Omega)$ and
$$
 \lim_{\eps\to0}\widetilde{\mathcal E}_\eps(u_\eps)= \lim_{\eps\to0}\mathcal E_\eps(u_\eps)= \surftens
\int_{\Omega\cap\partial^* E}[1+\vert\Hv_{\partial^* E}\vert^2]\,d\Ha^1+2\surftens L<+\infty.
$$
Therefore we can conclude that 
$$
\max\{\Gamma(L^1(\Omega))-\lim_{\eps\to 0}\widetilde{\mathcal E}_\eps(u),\,\Gamma(L^1(\Omega))-\lim_{\eps\to 0}\mathcal E_\eps(u)
\}<+\infty.
$$
For this choice of $u\in BV(\Omega,\{-1,1\})$ we expect that neither \eqref{lamb}, nor its analogue for $\Galim\Eeps(u)$ hold.  In fact: on the one hand \eqref{eq:RS-intro} and \eqref{eq:BeMu-intro} hold as soon as we localize the functionals $\widetilde{\mathcal E}_\eps$ and $\mathcal E_\eps$ on any open subset $\omega$ such that $\overline\omega\cap\{p_0,p_1\}=\emptyset$; on the other hand we cannot have 
$$
\lim_{\eps\to 0}\mu_\eps=\Ha^{1}_{\res \partial^* E} \text{ as Radon measures on }\Omega,
$$ 
 as this would contradict the fact (established in \cite{RoSch, ToneYuko, BeMu:09} and recalled above) that the rectifiable varifold associated with  the limit of the $\mu_\eps$ has $L^2$-bounded first variation.  Hence we expect that for every sequence $\{u_\eps\}_\eps\subset C^2(\Omega)$ such that $u_\eps\to u$ in $L^1(\Omega)$ we have
$$
\frac{1}{\surftens}\min\{\liminf_{\eps\to 0}\widetilde{\mathcal E}_\eps(u_\eps),\,\liminf_{\eps\to 0}\mathcal E_\eps(u_\eps)\}\geq\int[1+\vert\Hv_\mu\vert^2]\,d\mu>\int_{\partial^*E}[1+\vert\Hv_{\partial^* E}\vert^2]\,d\Ha^1,
$$
that is the last term on the right hand side is a too rough (or ``non-optimal'') lower-bound for  both $\widetilde{\mathcal E_\eps}(u_\eps)$ and $\mathcal E_\eps(u_\eps)$.
It is thus rather natural to try to answer the question: what are the $\Gamma$-limits of $\{\widetilde{\mathcal E}_\eps\}_\eps$ and $\{\mathcal E_\eps\}_\eps$ out of ``smooth sets''? 

We try to answer this question in the case $d=2$ only, and from now on, throughout the paper we will assume that $d=2$, unless otherwise specified. 

Since $\Gamma$-limits are necessarily lower semi-continuous functionals (see \cite[Proposition 4.16]{DM}), in view of \eqref{eq:RS-intro}, \eqref{eq:BeMu-intro} and \eqref{eq:2d}, a natural candidate for the $\Gamma$-limit of both $\widetilde{\mathcal E_\eps}$ and $\mathcal E_\eps$ is  the lower semi-continuous envelope (with respect to the $L^1(\Omega)$-topology) of the functional
\begin{gather}\label{eq:def-ela}
\mathcal F_o:L^1(\Omega)\to\,[0,+\infty],
\notag
\quad
u\mapsto \begin{cases} \int_{\Omega\cap\partial E}[1+\vert\Hv_{\partial E}\vert^2]\,d\Ha^1 & \text{if }u=\ind_E \text{ and }\Omega\cap\partial E\in C^2, 
\\
+\infty & \text{otherwise on } L^1(\Omega),
\end{cases}
\end{gather}
 that is the functional
\begin{align*}
\overline{\mathcal F}_o(u):=&\inf\{\liminf_{k\to\infty}\mathcal F_o(u_k): ~ L^1(\Omega)-\lim_{k\to\infty} u_k=u\}
\\
=&\sup\{\mathcal G(u):~ \mathcal G\leq\mathcal F_o \text{ on }L^1(\Omega),\, \mathcal G\text{ is lower semi-continuous on } L^1(\Omega)\}.
\end{align*}
Since by \cite[Theorem 3.2]{BDMP}  we have $\overline{\mathcal F}_o(u)=\mathcal F_o(u)$ for every $u=\ind_E$ such that $\Omega\cap\partial E\in W^{2,2}$ (see also \cite{Sch:09} for a more general statement),  by  \eqref{eq:RS-intro}, \eqref{eq:BeMu-intro} and the definition of $\overline{\mathcal F}_o$, we can conclude that
\begin{equation*}
\Gamma(L^1(\Omega))-
\lim_{\eps\to 0}\widetilde{\mathcal E}_\eps\leq\surftens\overline{\mathcal F}_o\text{ on }L^1(\Omega),
\quad
\Gamma(L^1(\Omega))-
\lim_{\eps\to 0}\mathcal{E}_\eps\leq\surftens\overline{\mathcal F}_o\text{ on }L^1(\Omega).
\end{equation*}
 We can now rephrase the results we obtain as follows:  we prove that $\Galim\mathcal E_\eps=\surftens\overline{\mathcal F}_o$ (at least under suitable boundary conditions for the phase-field variable), and we show that in general $\Galim\Etildeps <\surftens\overline{\mathcal F}_o$.

The outline of the paper is the following.
In Theorem \ref{theo:approx}, we show that the assumption
$$
\sup_{\eps>0}\mathcal E_\eps(u_\eps)<+\infty,
$$
implies additional ``regularity'' on the support of the measure $\mu:=\theta\Ha^1_{\res M}$ arising as  limit of the energy density measures $\mu_\eps$ defined above. 
Namely we establish that  in every point of $M\cap\Omega$ a (unique) tangent-line to $M$ is well defined. Moreover, in Corollary \ref{cor:sheldon-brown}, we show that  if $\{u_\eps\}_\eps\subset X$
where
$$
X:=\{u\in C^2(\Omega):~u(x)\equiv 1,\,\partial_{\nu_\Omega} u(x)=0,\,\forall x\in\partial\Omega\},
$$
then the set $M\cup\partial\Omega$ has an uniquely defined tangent line in every point.  
In view of  \cite{BeMu:07, BeMu:09} (see also  Theorem \ref{theo:rel-el} and Theorem \ref{the:mio} in this paper) this allows us to conclude that 
 $\Gamma(L^1(\Omega))-\lim_{\eps\to 0}\mathcal E_{\eps\res X}=\surftens\overline{\mathcal F}$, where $\overline{\mathcal F}$ denotes the $L^1(\Omega)$-lower semi-continuos envelope of  $\mathcal F_{o\res\mathscr K}$ and 
 $$
 \mathscr{K}:=\{E\subset\Omega:~E\text{ is open, compactly contained in }\Omega\text{ and }\partial E\in C^2\}
 $$
 (the restriction of $\mathcal F_o$ to $\mathscr K$ is a consequence of the fact that here we are considering the $\Gamma$-limit of $\mathcal E_{\eps\res X}$).
 
\noindent We remark that we do not expect that an analogoue of Theorem \ref{theo:approx} holds in space dimensions $d>2$.  In fact, to prove Theorem \ref{theo:approx} we make use of a blow-up argument and of some regularity results obtained in  \cite{Hutch-reg:84}, that are valid only for generalized $(d-1)$-dimensional hypersurfaces (namely, Hutchinson's curvature varifolds) with $p$-integrable (generalized) second fundamental form for some $p>(d-1)$.  Moreover, though we expect that an analogue of Corollary \ref{cor:sheldon-brown} holds also (at least) when $d=3$, to prove such a result we would probably need a different approach. In fact, in the proof of Corollary \ref{cor:sheldon-brown} we make an essential use of an ``explicit'' representation of $\overline{\mathcal F}$, that has been  established in \cite{BeMu:07} and is available only in two space dimensions.
 
 For what concerns the sequence  $\{\Etildeps\}_\eps$, in Theorem \ref{theo:sega} we show that in general the support of the limit measure does not necessarily have an unique tangent line in \textit{every} point, and we obtain the existence of a function $u=\ind_E\in BV(\Omega,\{-1,1\})$ such that 
$$
\Gamma(L^1(\Omega))-\lim_{\eps\to 0}\widetilde{\mathcal E}_\eps(u)<\surftens\overline{\mathcal F}_0(u)=\Galim\Eeps(u)=+\infty.
$$
This means that  the sharp interface limit of $\Eeps$ and $\Etildeps$ do not coincide as functionals on $L^1(\Omega)$, although as we already remarked $\Galim\Etildeps(u)=\Galim\Eeps(u)$ whenever $u=\ind_E$ and $\Omega\cap\partial E\in C^2$.

Although we are not able to completely identify the $\Galim\Etildeps(u)$ 
we shortly discuss how the results of \cite{DelPino} can be applied to obtain the value of the $\Gamma$-limit in a certain number of cases.

The paper is organized as follows. In Section \ref{labebysitter} we fix some notation, and recall some results about varifolds and the  lower semi-continuous envelope of $\mathcal F_o$. In Section \ref{sec:none}, for the readers convenience, we state the main results of \cite{RoSch, ToneYuko, BeMu:09}. In Section \ref{sec:puppamelo} we state our main results, the proofs of which are presented in Sections \ref{appotta}-\ref{smoothie}.  Finally, in Section \ref{baffing} we collect some preliminary lemmata needed in the proof of our main results.
  
\subsection*{Acknowledgment}
I wish to thank Giovanni Bellettini, Matthias R\"oger and Andreas R\"atz for several interesting discussions on the subject of this paper.
\section{Notation and preliminary results}\label{labebysitter}%

\subsection{General Notation}\label{iltroio}

Throughout the paper we adopt the following notation.
By $\Omega$ we denote an open bounded subset of $\R^2$ with smooth boundary. 
By $B_R(x):=\{z\in\R^2:~\vert z\vert<R\}$ we denote the euclidean open ball of radius $R$ centered in $x$. 

By $\mathcal L^2$ we denote the $2$-dimensional Lebesgue-measure, and by $\Ha^1$ the one-dimensional Hausdorff measure.

For every  set $E\subseteq \R^2$ we denote by $\chi_E$ the
characteristic function of $E$, that is $\chi_E(x)=1$ if $x \in E$,
$\chi_E(x)=0$ if $x \notin E$. Moreover we define the function $\ind_E$ by $\ind_E:=2\chi_E-1$. We denote by  $\overline
E$ and $\partial E$ respectively  the closure and
the topological boundary of $E$. All sets we consider are assumed
to belong to $\mathcal M$, the class of all measurable subsets of
$\R^2$.

We say that $E \subset \R^2$ is of class $W^{2,2}$
(resp. $C^k$, $k \geq 1$) in $\Omega$, and write $\Omega\cap\partial E\in W^{2,2}$ (resp. $\Omega\cap\partial E\in C^k$) if $E\cap\Omega$ is open and
can be locally represented as the subgraph of a
function of class $W^{2,2}$ (resp. $C^k$).

We say that a set $E\subset\R^2$ has finite perimeter in $\Omega$ if $\chi_E\in BV(\Omega)$, moreover if $E$ has finite perimeter by $\partial^*E$ we denote its reduced boundary (see \cite{Sim}).

We endow the space of  the $(2\times 2)$ 
matrices 
$M=(m_{ij})\in\R^{2\times 2}$
(resp. $2^3$ tensors $T=(t_{ijk})\in\R^{2^3}$) with 
the norm
\begin{equation}\label{amantedellosplit}
\vert M\vert^2 := {\rm tr}(M^T M) =\sum_{i,j=1}^2 (m_{ij})^2 \qquad 
\left({\rm 
resp.} 
~\vert T\vert^2:=\sum_{i,j,k=1}^2(t_{ijk})^2 \right),
\end{equation}
where $M^T$ is the transposed of $M$. 

If $P\in\R^{2\times 2}$ is a (symmetric) orthogonal projection matrix onto some 
subspace of $\R^2$ and $M$ is 
symmetric, then
\begin{equation}\label{eq:Kill-Bill}
\vert P^T M P\vert^2\leq\vert M\vert^2.
\end{equation}
%
\subsection{Differential Geometry}\label{prel:diff-geom}
Let $\Sigma$ be a smooth, compact oriented curve without boundary 
embedded in $\R^2$. 
If $x\in\Sigma$, we denote by $P_\Sigma(x)$ the orthogonal
projection onto the tangent line $T_x\Sigma$ to $\Sigma$ at $x$. 
Often we identify the linear operator $P_\Sigma(x)$  with the symmetric 
$(2\times 2)$-matrix ${\rm Id}-\nu_x\otimes\nu_x$ where $x \to \nu_x \in 
(T_x\Sigma)^\perp$ is a smooth unit 
covector field 
orthogonal to $T_x\Sigma$.

Let us recall that, when $\Sigma$ is given as a level surface $\{v=t\}$ 
of a smooth function $v$ such that $\nabla v\neq 0$ on $\{v=t\}$,  we 
can take at $x \in \{v=t\}$
\begin{equation*}
\nu_x=\frac{\nabla v(x)}{\vert\nabla v(x)\vert},\qquad P_\Sigma(x)
={\rm Id}-\frac{\nabla v(x)\otimes\nabla v(x)}{\vert\nabla v(x)\vert^2}.
\end{equation*}
The second fundamental form ${\bf B}_\Sigma$ of $\Sigma$ has the 
expression
\begin{equation*}
\sff_{\Sigma}=
\Big( P_\Sigma^T 
\frac{\nabla^2v}{\vert\nabla v\vert}
P_\Sigma\Big)
\otimes \frac{\nabla v}{\vert\nabla v\vert},
\end{equation*}
where $P_\Sigma^T = (P_\Sigma)^T$.
The definition of $\sff_\Sigma$ depends only on $\Sigma$ and not 
on the particular choice of the function $v$. Moreover 
$\sff_\Sigma(x)$,
if restricted to $T_x \Sigma$ and 
considered as a bilinear map from $T_x \Sigma \times T_x \Sigma$ 
with values in $(T_x \Sigma)^\perp$, 
 coincides with the usual notion of second fundamental form.
By 
$$
\Hv_\Sigma(x)=(H_1(x),H_2(x)) = {\rm tr}\Big( P_\Sigma^T 
\frac{\nabla^2v}{\vert\nabla v\vert}
P_\Sigma\Big) \nu_x,
$$
we denote the curvature vector of $\Sigma$ at $x $. 

Let us also define $A^\Sigma(x):=(A^\Sigma_{ijk}(x))_{1\leq i,j,k\leq 
3}\in\R^{2^3}$ as
\begin{equation}\label{eq:def-Aijk-Sigma}
A^\Sigma_{ijk}=\delta^\Sigma_i P_{\Sigma jk} \qquad {\rm on} ~ \Sigma,
\end{equation}
where $\delta_i^\Sigma := P_{\Sigma ij} \frac{\partial}{\partial x_j}$.

To better understand definition \eqref{eq:def-Aijk-Sigma}, 
it is useful to recall the 
links between 
$\sff_\Sigma$ and $A^\Sigma$ (see \cite[Proposition 5.2.1]{Hu}).

\begin{proposition}\label{prop:A-vs-B}
Set $A = A^\Sigma$, ${\bf B} = {\bf B}_\Sigma$
and ${\bf H} = {\bf H}_\Sigma$.
For $i,j,k \in \{1,2\}$
the following relations 
hold:
\begin{align}
& B^k_{ij}=P_{jl}A_{ikl},
\label{Bjork}
\\
& A_{ijk}=B^k_{ij}+B^j_{ik},
\label{eq:AvsB}
\\
& \Hv_i=A_{jij}=B^j_{ji}+B^i_{jj}.
\end{align}
\end{proposition}

Let $u \in C^2(\Omega)$.  We will often look at geometric 
properties
of the \textit{ensemble of the level sets} of $u$. 
We 
define
\begin{equation}\label{eq:def-proj-nablaperp}
\niu:=\frac{\nabla u}{\vert\nabla u\vert},\qquad 
\proju:={\rm Id}-\niu \otimes\niu, \qquad
\text{on }\{\grad u\neq 0\} ,
\end{equation}
 and  $\niu:=\mathbf{e}_3,\, \proju:={\rm 
Id}-\mathbf{e}_3\otimes\mathbf{e}_3$ on $\{\grad u=0\}$. Moreover we 
define 
the second fundamental form of the ensemble of the level sets of $u$ 
by
\begin{equation}\label{eq:def-sff-lev-surf}
\sff_u=\frac{(\proju)^T\nabla^2 u \proju}{\vert\nabla u\vert}\otimes\niu,
\end{equation}
on $\{\grad u\neq 0\}$ 
and $\sff_u:=\otimes^3\mathbf{e}_3$ on $\{\grad u=0\}$. 
Similarly we define
\begin{equation}\label{def:Aijk-u}
A_{ijk}^u:=-\proju_{il}\big[\partial_{l}((\niu)_j(\niu)_k)\big],
\end{equation}
on $\{\nabla u\neq 0\}$ and 
$A^u:=\otimes^3\mathbf{e}_3$ on $\{\nabla u=0\}$.

It will be convenient to consider $\sffu$ and $A^u$ 
as maps  defined on 
$G_2(\Omega)$ by 
$\sff_u(x,S):=\sff_u(x)$,
$A^u(x,S):=A^u(x)$.

\subsection{Geometric Measure Theory: varifolds}\label{luponeassassino}
Let us recall some basic fact in the theory of 
varifolds, the main bibliographic sources being \cite{Sim} and \cite{Hu}.

Let $M\subset\R^2$ be a Borel-set. We say that $M$ is $1$-rectifiable if there exists a countable family of  graphs (suitably rotated and translated) $\{\Gamma_n\}_{n\in\N}$ of Lipschitz functions $f_n$ of one variable such that $\Ha^1(M\setminus \cup_{n\in\N}\Gamma_n)=0$ and $\Ha^1(M)<+\infty$. 

By $G_{1,2}$ we denote the Grasmannian of $1$-subspaces of $\R^2$. We identify $T\in G_{1,2}$ with the projection matrix  $P_T\in\R^{2\times 2}$ on $T$, and endow $G_{1,2}$ with the relative distance as a compact subset of $\R^{2\times 2}$. Moreover, given $\Omega\subset\R^2$ open, we define the product space $G_1(\Omega):=\Omega\times G_{1,2}$, and endow it with the product distance.
    
We call \textit{varifold} any positive 
Radon measure on $G_1(\Omega)$. 
In this paper we are confined to curves, hence 
we use the terms varifold  to mean a $1$-varifold in $\Omega$.

By \textit{varifold convergence} we
mean the convergence as Radon measures
 on $G_1(\Omega)$.

For any varifold $V$ we define 
$\mu_V$ to be the Radon  measure on $\Omega$  
obtained by projecting $V$ onto $\Omega$.

Let $\rectifiableset
$ be a $1$-rectifiable subset of $\R^2$ and let  $\theta:\rectifiableset\to\R^+$ be a 
$\Ha^1 \res \rectifiableset$-measurable 
functions. We define the \textit{rectifiable varifold}
$
V=\var(\rectifiableset
,\theta),
$
by
\begin{gather*}
V(\phi):=\int_{\rectifiableset
}
\phi(x,T_x \rectifiableset
)\,\theta(x)d\Ha^2 \qquad\forall\phi\in 
C^0_c(G_2(\Omega)).
\end{gather*}
When $\theta$  takes values in $\N$
 we say that $V=\var(\rectifiableset
,\theta)$ is a  \textit{rectifiable integral varifold}
and we write $V \in \mathbf{IV}_1(\Omega)$.

Let $V$ be a
varifold on $\Omega$. We define \textit{the first variation 
of} $V$ as the linear operator
\begin{equation*}
\delta V:C^1_c(\Omega,\R^3)\to\R,\qquad Y\to \int {\rm tr}(S\nabla 
Y(x))\,dV(x,S).
\end{equation*}

We say that $V$ has \textit{bounded first variation}  
 if $\delta V$ 
can be 
extended to a linear continuous operator on  $C^0_c(\Omega,\R^2)$. In this case by $\vert \delta V\vert$ we 
denote the total variation of $\delta V$.
Whenever the 
varifold $V$ has bounded first variation we  call \textit{generalized 
mean curvature vector of} $V$ the vector field 
$$
\Hv_{V}=\frac{d\delta 
V}{d\mu_V},
$$
where the right-hand side denotes the Radon-Nikodym derivative of $\delta V$ with respect to $\mu_V$. 

\begin{remark}\rm\label{rem:local}
Let us recall  that if $V\in\mathbf{IV}_1(\Omega)$ and $V$ has bounded first variation then, by the results recently proved in \cite{Men, Masnou}, we have: the support of $\mu_V$ is a $1$-dimensional $C^2$-rectifiable subset of $\Omega$; $\Hv_V$ depends only on the local structure of the  varifold $V$, that is for every $V_1,\,V_2\in\mathbf{IV}_1(\Omega)$ we have $\Hv_{V_1}(x) =\Hv_{V_2}(x)$ for $\Ha^1$-a.e. $x\in \mathrm{spt}(\mu_{V_1})\cap \mathrm{spt}(\mu_{V_2})$.
\end{remark}

We say that a varifold $V$ is stationary if $\delta V\equiv 0$.

We say that $V\in\mathbf{IV}_1(\Omega)$ has $L^p$-\textit{bounded first variation} ($p>1$) if 
$$
\sup_{\substack{Y\in C^1_c(\Omega),\\ \|Y\|_{L^p(\mu_V)}\leq 1}}
\delta V(Y)<+\infty.
$$
It can be easily checked that every  $V\in\mathbf{IV}_1(\Omega)$ with $L^p$-bounded first variation verifies $\vert\delta V\vert<<\mu_V$ (as Radon measures), so that 
$$
\delta V(Y)=\int\Hv_V\cdot Y\,d\mu_V,\quad\Hv_V\in L^p(\mu_V).
$$

For every $V\in\mathbf{IV}_1(\Omega)$ with $L^p$-bounded first variation we set 
\begin{gather*}
\mathcal F_p(V):=\int[1+\vert\Hv_V\vert^p]\,d\mu_V
=\mu_V(\Omega)+\Big(\sup_{\substack{Y\in C^1_c(\Omega),\\ \|Y\|_{L^p(\mu_V)|}\leq 1}}
\delta V(Y)\Big)^p.
\end{gather*}

\begin{remark}\label{rem:briosch}\rm
If $V\in\mathbf{IV}_1(\Omega)$ has $L^p$-bounded first variation for some $p>1$, by \cite[Corollary 17.8]{Sim},  the $1$-density of $\mu_V$ in $x$ 
\begin{gather}\label{eq:den:bound}
\Theta(\mu_V,x):=\lim_{\rho\to 0}\frac{\mu_V(B_\rho(x))}{\pi\rho},
\end{gather}
is well defined everywhere on $\mathrm{spt}(\mu_V)$,  $\Theta(\mu_V,x)\in\N$ and $\Theta(\mu_V,x)<C$, where  $C>0$ is a constant that depends only on $\|\Hv_{\mu_V}\|_{L^p(\mu_V)}$. Moreover we can write $V=v(M,\theta)$ where $M=\mathrm{spt}(\mu_V)\cap\Omega$ and $\theta(x)=\Theta(\mu_V,x)$. In the rest of the paper we will always assume that varifolds with $L^p$-bounded first variation are represented in this manner. Eventually let us also recall that for $\Ha^1$-almost every $x_0\in \mathrm{spt}(\mu_V)$, there exists $P\in G_{1,2}$ such that
$$
\lim_{\rho\to 0^+}\frac{1}{\rho}\int\phi(\rho x+x_0,S)\,dV(x,S)=\theta(x_0)\int_P\phi(y,P)\,d\Ha^1,\quad\forall\phi\in C^0_c(G_1(\Omega).
$$
Moreover $P$ is a classical tangent line to $M$ at $x_0$ in the sense that
$$
\limsup_{\rho\to 0}\left\{\frac{\dist(x,P+x_0)}{\rho}:~x\in M\cap B_\rho(x_0)\right\}=0.
$$ 
\end{remark}

For our purposes we also need to introduce a further class of varifolds.  
Following \cite{Hu} we define the notion of Hutchinson's
curvature varifold with generalized second fundamental form.

\begin{definition}\label{def:Hutch-var}
Let $V \in {\bf IV}_1(\Omega)$.
We say that 
$V$ is a curvature varifold with generalized second 
fundamental form in  $L^p$ ($p>1$), if there exists $A_V=A^V_{ijk}\in 
L^p(V,\R^{2^2})$  such that for every 
function $\phi\in C_c^1(G_1(\Omega))$ and $i=1,2$, 
\begin{equation}\label{eq:def-Hutch-var} 
\int_{G_{2}(\Omega)}(S_{ij}\partial_j 
\phi+A^V_{ijk}D_{m_{jk}}\phi+A^V_{jij}\phi)\, 
dV(x,S)=0,
\end{equation}
where $D_{m^{}_{jk}} \phi$ denotes the derivative of $\phi (x,\cdot)$ 
with respect to its $jk$-entry variable.

Moreover we define the generalized second fundamental form $\sff_V
=(B^k_{ij})_{1\leq i,j,k\leq 3}$ of $V$ as
\begin{equation}\label{eq:def-gen-sff}
B^k_{ij}(x,S):=S_{jl}A^V_{ikl}(x,S).
\end{equation}
Eventually  by $\mathscr{CV}^p_1(\Omega)$ we denote the calss of Hutchinson's curvature varifolds with $p$-integrable second fundamental form in $\Omega$. 
\end{definition}
\begin{remark}\rm
If $V=\var(\Sigma,1)$, where $\Sigma$ is a smooth, compact 
surface without boundary, the generalized second fundamental form as well as the mean curvature and the tensor $A_V$ coincide with the classical quantities defined in Section \ref{prel:diff-geom}. Moreover, for every $V\in \mathscr{CV}^p_1(\Omega)$ ($p>1$)  the functions $A^V_{ijk},\, B^k_{ij}$ verify $V$-a.e.the identities stated in Proposition \ref{prop:A-vs-B}. 
\end{remark} 
\begin{remark}\label{rem:MC-Hutch-MC-Allard}\rm
Every curvature varifold $V$ with generalized  second
fundamental form in $L^p$ has also $L^p$-bounded first variation and
\begin{equation}\label{eq:MC-Hutch-MC-Allard} 
\Hv_V(x)=(A_{212}(x,T_x\mu_V),A_{121}(x,T_x\mu_V))\in L^p(\mu_V,\R^2),
\end{equation}
for $\mu_V$ almost every $x \in \Omega$ (see \cite{Hu}).  Moreover if $V\in\mathscr{CV}^p_1(\Omega)$, by Proposition \ref{prop:A-vs-B}, we have
\begin{gather}
\mathcal F_p(V)=\int[1+\vert\Hv_V\vert^p]\,d\mu_V=
\int[1+\vert\sff_V\vert^p]\,dV=\int[1+\vert A_V\vert^p]\,dV.
\label{eq:astianatte}
\end{gather}
Let us also recall that there are, however, varifolds $V
\in\mathbf{IV}_1(\Omega)$ with $L^p$-bounded first variation for every $p>1$, that do not belong to $\mathscr{CV}^p_1(\Omega)$. An example is given by the (stationary) varifold $\var(M,1)\in B_R$ where $M$ is given by the union of three line segments of length $R$ having one end point in the origin, and forming angles of $2\pi/3$ radiants one with the other. 
\end{remark}

Eventually we introduce the set $\mathscr{D}(\Omega)\subsetneq\mathscr{CV}^2_1(\Omega)$ of  Hutchinson's curvature varifolds that can be approximated (in the varifolds topology) by a sequence of $C^2$-smooth embedded curves in $\Omega$, having uniformly $L^2$-bounded second fundamental form. More precisely we give the following
\begin{definition}\label{def:D}
 We define the set $\mathscr D(\Omega)$ as the set of $V\in\mathscr{CV}^2_1(\Omega)$ for which there exists a sequence $\{E_k\}_k$ of open, bounded subsets  with smooth boundary such that $E_k\subset\subset\Omega$ and such that 
\begin{gather*}
\lim_{k\to\infty}\var(\partial E_k,1)=V,\text{ as varifolds in }\Omega,
\\
\sup_{k\in\N}\mathcal F_2(V_k)=\sup_{k\in\N}\int_{\partial E_k}[1+\vert\Hv_{\partial E_k}\vert^2]\,d\Ha^1=\sup_{k\in\N}\int_{\partial E_k}[1+\vert\sff_{\partial E_k}\vert^2]\,d\Ha^1<+\infty.
\end{gather*}
\end{definition} 
As a straightforward consequence of the results proved in \cite{BeMu:07} we have the following characterization 
\begin{equation}\label{eq:Cahr-D}
\mathscr D(\Omega)=\Big\{V=\var(M,\theta)\in\mathscr{CV}^2_1(\R^2):~ M\cup\partial\Omega\text{ has an unique tangent line in every point}\}.
\end{equation}

\begin{remark}\label{rem:up-to-bdry}\rm
If in Definition \ref{def:D} we drop the assumption $E_k\subset\subset\Omega$ on the sequence of smooth sets approximating $V=\var(M,\theta)$, then  \eqref{eq:Cahr-D} ceases to hold. In fact, in this case $M$ has a unique tangent line in every point belonging to $M\cap\Omega$ (see Proposition \ref{lem:labradford} below), but there might be points $p\in M\cap\partial \Omega$ where the tangent line to $M\cup\partial\Omega$ is not unique. As a consequence, though $V\in\mathscr{CV}^2_1(\Omega)$,  in general we have $V\notin\mathscr{CV}^2_1(\R^2)$. 
\end{remark}

We conclude this section with a further easy consequence of \cite{BeMu:07}, that we need in the proof of Theorem \ref{theo:sega}.

\begin{proposition}\label{lem:labradford}
Let $V=\var(M,\theta)$ be an integrable, rectifiable varifold with $L^2$-bounded first variation in $\Omega$. Suppose we can find a sequence of  manifolds $M_k$ smooth, embedded and without boundary in $\Omega$, such that   $V=\lim_{k\to\infty}\var(M_k,1)$ with respect to  varifolds convergence in $\Omega$ and such that
\begin{equation}\label{silvioingalera}
\sup_{k\in\N}\int_{M_k}1+ \vert\Hv_{M_k}\vert^2\,d\Ha^1<+\infty.
\end{equation}
Then $\mathrm{spt}(\mu_V)=M$ has an unique tangent line in every point of $M\cap U$ for every $U\subset\subset\Omega$.  
\end{proposition}

\begin{remark}\label{rem:maiala}\rm
Let us mention that we expect that the arguments used in \cite{BeMu:07} can be adapted to prove also the converse of Proposition \ref{lem:labradford}. That is, if $\var(M,\theta)\in\mathscr{CV}^2_1(\Omega)$ is such that $M$ has an unique tangent line in every point of $M\cap\Omega$ then there exists a sequence of  manifolds $M_k$ smooth, embedded and without boundary in $\Omega$, such that   $V=\lim_{k\to\infty}\var(M_k,1)$ with respect to  varifolds convergence in $\Omega$, and such that $V_k$ verify \eqref{silvioingalera}.
\end{remark}

\begin{remark}\label{rem:chloe}
Let $V=\var(M,\theta)\in\mathscr{CV}^2_1(\Omega)$. In \cite{BeMu:07} it has been proved that to say that in every point of $M\cap \Omega$ an unique tangent line is well defined, is equivalent to say that  $M\cap\Omega$ can be locally (and up to rigid motions) represented as a finite union of graphs of $W^{2,2}$-functions that do not cross each other.
\end{remark}  

\subsection{Preliminary Results on the Relaxed elastica Functional}
Let us define the functional
\begin{gather}
\mathcal F:=\mathcal F_{o\res \mathscr K}:L^1(\Omega)\to\,[0,+\infty],
\notag
\\
u\mapsto \begin{cases} \int_{\partial E}[1+\vert\Hv_{\partial E}\vert^2]\,d\Ha^1 & \text{if }u=\ind_E, \text{ and }E\subset\subset\Omega,\,\partial E\in C^2,  
\\
+\infty & \text{otherwise on } L^1(\Omega),
\end{cases}
\label{eq:def-el}
\end{gather}
and its $L^1$-lower-semicontinuous-envelope 
\begin{gather}
\overline{\mathcal F}(u):=\inf\{\liminf_{k\to\infty}\mathcal F(u_k):~\lim_{k\to\infty} u_k= u\text{ in }L^1(\Omega)\}.
\label{figa}
\end{gather}
\begin{remark}\rm\label{rem:no-confusion}
We remark that if $E\subset\Omega$ is open and with smooth boundary, $u:=\ind_E$  and $V=\var(\partial E,1)\in\mathbf{IV}_1(\Omega)$, we have $\mathcal F_2(V)=\mathcal F(u)$.
\end{remark} 

As a straightforward consequence of  \cite[Theorem 4.3]{BeMu:07} we have the following
\begin{theorem}\label{theo:rel-el}
Let $E\subset\Omega$ and $u=\ind_E\in L^\infty(\Omega,\{-1,1\})$. Then $\overline{\mathcal F}(u)<+\infty$ if and only if $u\in BV (\Omega,\{-1,1\})$ and the set
\begin{align*}
\mathscr A(E):=\{V=\var(M,\theta)\in\mathscr{D}(\Omega):&~M\supset\partial^*E\neq\emptyset,
\\
&\theta(x) \equiv 1\,\mathrm{mod}2, \forall x\in\partial^*E,
\\
&\theta(x) \equiv 0\,\mathrm{mod}2, \forall x\in\mathrm{spt}(\mu_V)\setminus\partial^*E\},
\end{align*}
is not empty. Moreover, if $\mathscr A(E)\neq\emptyset$, the following representation formula holds
\begin{gather*}
\overline{\mathcal F}(u)=\min_{V\in\mathscr A(E)}\mathcal F_2(V).
\end{gather*}
In particular if $\partial E$ is $W^{2,2}$-smooth in $\Omega$ and, if $\partial E\cap\partial \Omega\neq\emptyset$, $\partial E$ touches  $\partial\Omega$ tangentially, then
$$
\overline{\mathcal F}(u)=\mathcal F(u).
$$
\end{theorem}
\begin{remark}\rm\label{rem:nadiacassini}
If $\{E_k\}_{k\in\N}$ is a sequence of open smooth subsets of $\Omega$ (that do not necessarily verify $E_k\subset\subset\Omega$) such that $L^1(\Omega)-\lim_{k\to\infty}u_k=\ind_{E}$, where $E$ is a subset with smooth boundary, by  \cite{BDMP,Sch:09} we still can conclude that  
$$
\liminf_{k\to\infty}\int_{ \Omega\cap\partial E_k}1+\vert\Hv_{\partial E_k}\vert^2\,d\Ha^1\geq \int_{\Omega\cap\partial E}1+\vert\Hv_{\partial E}\vert^2\,d\Ha^1.
$$
\end{remark}
\subsection{$\Gamma$-convergence}\label{bomboloni}
Let $X$ be a topological space and $F_\eps:X\to[0,+\infty]$ a sequence of functionals on $X$. We say that $F_\eps$
  $\Gamma$-converge to the $\Gamma$-limit $F:X\to[0,+\infty]$ in $X$, and we write $\Gamma(X)-\lim_{\eps\to 0}F_\eps=F$,
 if the following two conditions hold:
\begin{itemize}
\item Lower bound inequality (or $\Gamma-\liminf$-inequality): For every sequence $\{x_\eps\}\subset X$ such that $\lim_{\eps\to 0}x_\eps=x$ in $X$,
$$
\liminf_{\eps\to 0}
F_\eps(x_\eps)\geq F(x).
$$
\item Upper bound inequality (or $\Gamma-\limsup$-inequality): For every $x\in X$, there exists a recovery sequence $\{x_\eps\}_\eps\subset X$ such that
\begin{gather*}
\lim_{\eps\to 0}x_\eps=x\text{ in }X,
\quad
\limsup_{\eps\to 0}F_\eps(x_\eps)\leq F(x).
\end{gather*}
\end{itemize}

\section{Preliminary known Results on Diffuse Interfaces Approximations of $\mathcal F$}\label{sec:none}
We begin this section specifying some further notation needed in the sequel.

\noindent We set
\begin{equation*}
W(r):=\frac{1}{4} (1-r^2)^{2}, \qquad r \in \R,
\end{equation*}
and
\begin{equation}\label{dito}
\surftens:=\int_{-1}^1\sqrt{2W(s)}\,ds.
\end{equation}
If $\gamma(s):=\tanh(s)$ we have
$\ddot \gamma=\frac{d}{ds}
(W(\gamma))$,
\\
\begin{gather*}
\int_\R\vert\dot\gamma\vert^2\,ds=\int_\R 2W(\gamma)\,ds
=\surftens,
\end{gather*}
and
\begin{equation}\label{eq:def-surf-tens}
\surftens=\min\left\{\int_\R\left(\frac{\vert \dot v\vert^2}{2}+
W(v)\right)\,ds:\, v\in H^1_{{\rm loc}}(\R),\,
\lim_{s\to\pm\infty}v(s)=\pm1\right\}.
\end{equation}

To every sequence $\{u_\eps\}_\eps\subset C^2(\Omega)$ we associate 
\begin{itemize}
\item the sequences of Radon measures
\begin{gather}
\mu_\eps:=\Big(\frac{\eps}{2}\vert\nabla u_\eps\vert^2+\frac{W(u_\eps)}{\eps}\Big)\mathcal L^2_{\res\Omega},\quad \widetilde\mu_\eps:=\eps\vert\nabla u_\eps\vert^2\mathcal L^2_{\res\Omega};
\end{gather}
\item the sequence of diffuse varifolds
\begin{gather}\label{eq:ridi-culo}
\Vepsueps(\phi):=\surftens^{-1}\int\phi(x,\projueps(x))~d\widetilde\mu_\eps(x), \quad \forall 
\phi\in C^0_c(G_1(\Omega)),
\end{gather}
where $P^{u_\eps}(x)$  denotes the projection on the tangent space to the level line of $u_\eps$ passing through $x$ (see \eqref{eq:def-proj-nablaperp}).
\end{itemize}

The next result has been proved in \cite{RoSch, ToneYuko}
\begin{theorem}\label{the:RS} Let $\{u_\eps\}\subset C^2(\Omega)$ be a sequence such that 
\begin{equation}\label{eq:ass-RS}
 \sup_{0<\eps}\widetilde{\mathcal E_\eps}(u_\eps)= \sup_{0<\eps}
\mathcal P_\eps(u_\eps)+\mathcal W_\eps(u_\eps)< +\infty. 
\end{equation}
There exists a subsequence (still denoted by $\{u_\eps\}$)
converging to $u=\ind_E$ in $L^1(\Omega)$, where $E$ is a finite
perimeter set. Moreover 
\begin{itemize} 
\item[(A)]
$\mu_{\eps}\rightharpoonup\mu$ as $\eps\to 0^+$ weakly$^*$ in $\Omega$ 
as Radon measures and $\mu$ 
verifies 
\begin{equation*} \mu \geq \surftens \Ha^1 \res\partial E.
\end{equation*} 
In addition 
\begin{equation}\label{eq:the-RS-van-discr}
\lim_{\eps\to 0^+}\int_\Omega\vert\xi_\eps\vert\,dx=0,
\end{equation} 
where 
$$
\xi_\eps:=\Big(\frac{\eps}{2}\vert\nabla u_\eps\vert^2-\frac{W(u_\eps)}{\eps}\Big),
$$ 
and hence 
\begin{equation}\label{eq:lim-tuttuguale}
\mu = \lim_{\eps\to 0^+}\muepsueps=\lim_{\eps\to 
0^+}
\tildemuepsueps
=\lim_{\eps\to 0^+}\frac{2W(u_\eps)}{\eps}\LL^2 \res \Omega {\it ~as~
Radon~measures.} 
\end{equation} 
\item[(B)] The sequence $\{\Vepsueps\}$ 
 converges in  the varifolds 
sense to an
integral-rectifiable varifold $V=\var(M,\theta)\in\mathbf{IV}_1(\Omega)$ with $L^2$-bounded first variation, and such that 
$\mu_V = \surftens^{-1}\mu$. Moreover the function $\theta$ assumes odd (respectively even) values on $\partial^*E$ (respectively $M\setminus\partial^*E$). 
\item[(C)] For any $Y\in C^1_c(\Omega;\Rn)$ we have 
\begin{equation}\label{eq:the-RS-conv-first-var} 
\surftens\lim_{\eps\to 0^+}\delta
\Vepsueps(Y) =\lim_{\eps\to 0^+}\int_{\Omega} \Big(\frac{W^\prime(u_\eps)}{\eps}-\eps\Delta u_\eps\Big)\nabla 
u_\eps\cdot Y\,dx=
-\int \Hv_{V}\cdot Y\,d\mu, 
\end{equation}
and 
\begin{equation}\label{eq:the-RS_lsc-var} 
\surftens\mathcal F_2(V)=\surftens\int_{\Omega}\vert
\Hv_{V}\vert^2\,d\mu_V \leq 
\liminf_{\eps\to
0^+}\frac{1}{\eps}\int_{\Omega}\left(\eps\Delta u_\eps
-\frac{W^\prime(u_\eps)}{\eps}\right)^2\,dx.
\end{equation}
\end{itemize} 
\end{theorem}

As a straightforward consequence of \eqref{eq:the-RS_lsc-var},  Remark \ref{rem:local}, and \cite{BePa:93} we obtain the following
\begin{corollary}\label{cor:conv-smt}
For every $u=\ind_E\in BV(\Omega,\{-1,1\})$ such that $\Omega\cap\partial E\in W^{2,2}$, we have
$$
\Gamma(L^1(\Omega))-\lim_{\eps\to 0}\widetilde{\mathcal E}_\eps(u)
=\surftens\int_{\Omega\cap\partial E}[1+\vert \Hv_{\partial E}\vert^2]\,d\Ha^1=\surftens\mathcal F(u).
$$
\end{corollary}

Next we recall some of the main results obtained in \cite{BeMu:09} concerning
the $\Gamma$-convergence of the sequence $\mathcal E_\eps:=\mathcal P_\eps+\mathcal B_\eps$.
\begin{theorem}\label{the:mio}
Let $\{u_\eps\}\subset C^2(\Omega)$  be such that 
\begin{gather}\label{eq:mio}
\sup_{\eps>0}\mathcal E_\eps(u_\eps):=\sup_{\eps>0}\mathcal P_\eps(u_\eps)+\mathcal B_\eps(u_\eps)<+\infty.
\end{gather}
Then there exists a subsequence (still denoted by $\{u_\eps\}$)
converging to $u=\ind_E$ in $L^1(\Omega)$, where $E$ is a finite
perimeter set. Moreover 
\begin{itemize} 
\item[(A1)]
$\mu_{\eps}\rightharpoonup\mu$ as $\eps\to 0^+$ weakly$^*$ in $\Omega$ 
as Radon measures and $\mu$ 
verifies 
\begin{equation*} \mu \geq \surftens \Ha^1 \res\partial E.
\end{equation*} 
In addition 
\begin{equation}\notag
\lim_{\eps\to 0^+}\nabla\xi_\eps\,\LL^2_{\res\Omega}=0 {\it ~as~ Radon~measures},\quad \lim_{\eps\to 0}\|\xi_\eps\|_{L^p(\Omega)},{\it{~for~every~}}1<p<2,
\end{equation} 
and \eqref{eq:lim-tuttuguale} holds.
\item[(B1)] The sequence $\{\Vepsueps\}$ 
 converges to a varifold $V=\var(M,\theta)\in\mathscr{CV}^2_1(\Omega)$, such that: 
$\mu_V = \surftens^{-1}\mu$; and such that the function $\theta$ assumes odd (respectively even) values on $\partial^*E$ (respectively $M\setminus\partial^*E$). 
\item[(C1)] Let $A^u_{ijk}$ ($i,j,k=1,2$) be as in \eqref{def:Aijk-u}. For every $\phi\in C^1_c(G_1(\Omega);\R^2)$ we have 
\begin{equation}\label{eq:the-BeMu-conv-first-var} 
\lim_{\eps\to 0^+}\int_{\Omega} A^{u_\eps}_{ijk}(x,S)\phi(x,S)\,d\Vepsueps=\int A^V_{ijk}(x,S)\phi(x,S)\,dV(x,S), 
\end{equation}
for every $i,j,k=1,2$. Moreover 
\begin{equation}\label{eq:the-BeMu_lsc-var} 
\surftens\mathcal F_2(V) =\surftens\int\vert \sff_V\vert^2\,dV\leq
\liminf_{\eps\to
0^+}\frac{1}{\eps}\int_{\Omega}\left\vert\eps\nabla^2 u_\eps
-\frac{W^\prime(u_\eps)}{\eps}\niueps\otimes\niueps\right\vert^2\,dx.
\end{equation}
\end{itemize} 
\end{theorem}
\begin{remark}\rm\label{rem:inutile}
We notice that, in view of \eqref{smanettando}, the main assumption of Theorem \ref{the:mio}, that is \eqref{eq:mio}, is stronger than the main assumption of Theorem \ref{the:RS}, that is \eqref{eq:ass-RS}. 
 However also the conclusions of Theorem \ref{the:mio} are stronger than those of Theorem \ref{the:RS}. In fact, in Theorem \ref{the:mio}-(A1) the convergence to zero of  the discrepancies $\xi_\ep$  is proved to hold with respect to a topology  that is stronger than the one with respect to which the vanishing of the discrepancies is obtained in Theorem \ref{the:RS}-(A). Moreover in Theorem \ref{the:mio}-(B1) the limit varifold $V$ belongs to the set $\mathscr{CV}_1^2(\Omega)$, which is  strictly contained in the set of varifolds $V\in\mathbf{IV}_1(\Omega)$ having $L^2$-bounded first variation (see Remark \ref{rem:MC-Hutch-MC-Allard}). 
\end{remark} 

Eventually we notice that, by a straightforward adaptation of the proof of Corollary \ref{cor:conv-smt} and \cite[Theorem 4.2]{BeMu:09},  we obtain
\begin{corollary}\label{cor:conv-smt-BeMu}
For every $u=\ind_E\in BV(\Omega,\{-1,1\})$ such that $\Omega\cap\partial E\in W^{2,2}$, we have
$$
\Gamma(L^1(\Omega))-\lim_{\eps\to 0}\mathcal E_\eps(u)
=\surftens\int_{\partial E}[1+\vert \sff_{\partial E}\vert^2]\,d\Ha^1=\surftens\int_{\partial E}[1+\vert \Hv_{\partial E}\vert^2]\,d\Ha^1=\surftens\mathcal F(u).
$$
\end{corollary}

\section{Main Results}\label{sec:puppamelo}

The first of our main results shows that  every varifold $V=\var(M,\theta)\in\mathscr{CV}_1^2(\Omega)$ arising as the limit of diffuse interface varifolds veryfing \eqref{eq:mio} (see Theorem \ref{the:mio}-(B1)) is more regular than a generic element of $\mathscr{CV}^2_1(\Omega)$. In fact we show that  $M$ has an unique tangent line at \textit{every} point $p\in M\cap\Omega$ (consequently $M$ can be  represented, locally and up to rigid motions,  as the finite union of the graphs of $W^{2,2}$-functions, see Remark \ref{rem:chloe}).
\begin{theorem}\label{theo:approx}
Let $\{u_\eps\}_\eps\subset C^2(\Omega)$ satisfy \eqref{eq:mio}. Let $\Vepsueps$ be as in \eqref{eq:ridi-culo} and suppose $\lim_{\eps\to 0}\Vepsueps=V=\var(M,\theta)\in\mathscr{CV}_1^2(\Omega)$. Then $M$ has an unique tangent line in every $p\in M\cap\Omega$.
\end{theorem} 

As a consequence of Theorem \ref{theo:approx} we obtain the following full $\Gamma(L^1)$-convergence result
\begin{corollary}\label{cor:sheldon-brown}
Let 
$$
X:=\{u\in C^2(\Omega):~u(x)\equiv 1,\,\partial_{\nu_\Omega} u(x)=0,\,\forall x\in\partial\Omega\}.
$$
Define (with a small abuse of notation)
\begin{gather*}
\mathcal E_{\eps\res X}:L^1(\Omega)\to\,[0,+\infty],
\quad
u\mapsto \begin{cases} \mathcal P_\eps(u)+\mathcal B_\eps(u) & \text{if }u\in X,
\\
+\infty & \text{otherwise on } L^1(\Omega)
\end{cases}.
\end{gather*}
 Then $\Gamma(L^1(\Omega))-\lim_{\eps\to 0}\mathcal E_{\eps\res X}=\overline{\mathcal F}$, where $\overline{\mathcal F}$ is as 
in \eqref{figa}.
\end{corollary}

\begin{remark}\rm\label{rem:bombamerd}
We remark that from the proof of Corollary \ref{cor:sheldon-brown} it follows that the $\Gamma$-limit of the sequence $\{\mathcal E_{\eps\res X}\}_\eps$ with respect to the varifold convergence of $\Vepsueps$ is given by the functional
\begin{equation*}
V\mapsto\begin{cases} \mathcal F_2(V) &\text{if }V\in\mathscr{D}(\Omega),
\\
+\infty &\text{otherwise on }\mathscr{CV}^2_1(\Omega),
\end{cases}
\end{equation*}
where $\mathcal F_2$ has been defined in \eqref{eq:astianatte}.
\end{remark}

\begin{remark}\rm\label{rem:strozzaminchie}
In Corollary \ref{cor:sheldon-brown} we need to introduce the space $X$ in order to constrain the ``diffuse interfaces'' 
$$
\Sigma_{\eps,\delta}:=\{x\in\Omega:~\vert u_\eps\vert <1-\delta\},
$$
to be compactly contained in $\Omega$ for every $\eps$ and $\delta$ positive. This fact, together with the results of Theorem \ref{the:mio}, enables us to conclude that the measure $\mu_V=\theta\Ha^1_{\res M}$ can be approximated by a sequence obtained restrcting the
 $\Ha^1$-measure to the boundaries of a sequence of open subsets compactly contained in $\Omega$, and in turn to apply Theorem \ref{theo:rel-el}. Proving a full $\Gamma$-convergence result when the functionals $\{\mathcal E_\eps\}_\eps$ are defined on a more general functions' space that allows the diffuse interfaces $\Sigma_{\eps,\delta}$ to hit the boundary, seems to be merely a technical point that can be solved by proving that the ``conjecture'' stated in Remark \ref{rem:maiala} is true.
\end{remark}

Eventualy we also obtain some results concerning the
$\ Gamma$-limit of the sequence of functionals $\{\widetilde{\mathcal E}_\eps\}_\eps$, and its relation with  $\overline{\mathcal F}$. More precisely in Theorem \ref{theo:sega}, as a quite direct consequence of the results proved in \cite{Dang-Fife-Pel} (see also \cite{{CabTe:JEMS, DelPino}}), we  prove the following
\begin{theorem}\label{theo:sega}
Let $\Omega\subset\R^2$. There exists a sequence $\{u_\eps\}_\eps\subset C^2(\Omega)$ such that  
\begin{gather*}
L^1(\Omega)-\lim_{\eps\to 0}u_\eps=u=\ind_E\in BV(\Omega,\{-1,1\})\text{ for some }E\neq\emptyset;
\\
 \lim_{\eps\to 0}\Vepsueps=V=\var(\partial E,1)\in\mathscr{CV}_1^2(\Omega), 
 \\
\sup_{\eps>0}\mathcal P_\eps(u_\eps)<+\infty, \quad \mathcal W_\eps(u_\eps)\equiv 0,
\end{gather*}
and  such that $\partial E=\mathrm{spt}(\mu_V)$ does not have an unique tangent line in every point.

Moreover 
\begin{gather}
\lim_{\eps\to 0}\widetilde{\mathcal E}_\eps(u_\eps)=\mathcal F_2(V)<\overline{\mathcal F}_0(u)=\Galim\Eeps(u)=+\infty.
\label{doggy}
\end{gather}
\end{theorem}
\begin{remark}\label{rem:moanin}
Although we are not able to identify the $\Gamma$-limit of the sequence $\{\Etildeps\}_\eps$, we believe that for any given varifold $V\in\mathscr{CV}^2_1(\Omega)$, combining the results of \cite{Hutch-reg:84} with those of \cite{DelPino} and \cite{BePa:93},  it is possible to construct (with some additional work) a sequence $\{u_\eps\}_\eps\subset C^2(\Omega)$ such that 
$$
\lim_{\eps\to 0}\Etildeps(u_\eps)=\mathcal F_2(V).
$$
\end{remark}

%
\section{Preliminary Lemmata}\label{baffing}

In order to prove Theorem \ref{theo:approx} we need the following Lemmata.

\begin{lemma}\label{lem:app-smooth-flat}
Let $\{V_k:=\var(M_k,1)\}_k\subset\mathbf{IV}_1(B_{2R})$. Suppose that $M_k\cap B_{2R}$ are  smooth $C^2$-embedded
$1$-manifolds without boundary in $B_{2R}$, and
\begin{equation}\label{eq:strong-flat}
\begin{split}
& 0<\liminf_{k\to\infty}\mu_{V_k}(B_{R})\leq\limsup_{k\to\infty}\mu_{V_k}(B_{2R})=K<+\infty,  
\\
& \lim_{k\to\infty}
\vert\delta V_k\vert(B_{2R})=\lim_{k\to\infty}\int_{M_k}\vert\Hv_{M_k}\vert\, d\Ha^1=0.
\end{split}
\end{equation}
There exist a finite collection of $1$-dimensional affine subspaces $T_1,\dots,T_N$ of $\R^2$ such that 
\begin{gather}
T_i\cap T_j\cap B_R=\emptyset,\quad \text{for }i\neq j,\,i,j\in\{1,\dots,N\},
\label{eq:sheldon}
\end{gather}
 and a subsequence (not relabelled) $\{V_k\}_k\subset\mathbf{IV}_1(B_{2R})$ such that 
\begin{gather}
\lim_{k\to\infty}V_k(\phi)=\sum_{j=1}^N\Theta_j\int_{T_j}\phi(x,T_j)\,d\Ha^1=:V(\psi),
\quad \forall \psi\in C^0_c(G_1(B_R)),
\label{eq:cooper}
\end{gather}
where $\Theta_j\in \N$ are constants.
\end{lemma}
\proof
By \eqref{eq:strong-flat} we can apply Allard's compactness Theorem (see \cite[Theorem 42.7]{Sim}), and extract a subsequence such that $V_k\to V$, where $V\in \mathbf{IV}_1(B_{2R})$ is stationary  in $B_{2R}$, and $\mu_V(B_R)>0$.

Next we claim that (up to subsequences):
\begin{itemize}
\item[(i)] there are no closed curves between the connected components of $M_k\cap B_{2R}$;  
\item[(ii)] the connected components $M_k\cap B_{2R}$ intersecting $B_{3R/2}$ are in a fixed number.
\end{itemize}
 In fact, suppose that along a subsequence $\{M_{k^\prime}\}_{k^\prime}$ we  can find a closed curve $\widetilde M_{k^\prime}\subset M_{k^\prime}$ such that 
$\widetilde M_{k^\prime}\subset\subset B_R$  for every $k^\prime\in\N$.  Then
\begin{align*}
\vert\delta V_{k^\prime}\vert(B_{2R})=&\int_{M_{k^\prime}\cap B_{2R}}\vert\Hv_{M_{k^\prime}}\vert\,d\Ha^1 
\geq \int_{\widetilde M_{k^\prime}}\vert
\Hv_{\widetilde M_{k^\prime}}\vert\,d\Ha^1
\\
\geq& \Big\vert\int_{\widetilde M_{k^\prime}}\Hv_{\widetilde M_{k^\prime}}\,d\Ha^1\Big\vert=2\pi,
\end{align*}
which is in contradiction with \eqref{eq:strong-flat}.  Hence (i) holds.

Let us now prove (ii). Any $C^2$-embedded, non-closed curve without boundary in $B_{2R}$ intersecting $B_{3R/2}$ has a length of at least $R/2$. Hence the number of connected components of $M_k$ such that $M_k\cap B_R\neq \emptyset $ is smaller or equal than $2K/R$. Therefore, possibly passing to a further subsequence, we can suppose that the number of connected components of $M_k\cap B_{3R/2}$ equals a certain $N\in\N$ for every $k\in \N$.

In view of the results establised above and the assumption \eqref{eq:strong-flat}, we can find a constant $C>0$, a collection of $N$ intervals   $I_j\subset \R$ $j=1,\dots,N$ and $N$ sequence of maps $\{\alpha_{j,k}\}_{k\in\N}\subset C^2(I_j,B_{3R/2})$ such that, for $j=1,\dots,N$ and $k\in\N$, we have 
\begin{gather*}
C<\vert\dot{\alpha}_{j,k}\vert=const.\text{ on }I_j,
\quad
M_k\cap B_R=\bigcup_{j=1}^N(\alpha_{j,k})(I_j)\cap B_{R}.
\end{gather*}
Since
$$
\lim_{k\to\infty}\vert\delta V_k(B_R)\vert=\lim_{k\to\infty}\sum_{j=1}^N\frac{1}{l(\alpha_{j,k})}\int\vert\ddot\alpha_{j,k}\vert\,dt=0,
$$
  we have (up to the extraction of a further subsequence) $\alpha_{j,k}\to\alpha_j$ strongly in $W^{2,1}(I_j)$, for every $j=1,\dots,N$. Moreover by
 $$
 \sup_{s,t\in I_j}\vert\dot\alpha_{j,k}(s)-\dot\alpha_{j,k}(t)\vert\leq \sup_{s,t\in I_j}\int_s^t
\vert\ddot\alpha_{j,k}(\tau)\vert\,d\tau\leq\int_{I_j}\vert\ddot\alpha_{j,k}(\tau)\vert\,d\tau\to 0,
$$ 
we also have (again up to a subsequence) $\alpha_{j,k}\to\alpha_j$  uniformly and $\ddot{\alpha}_j\equiv0$ on $I_j$. Therefore $\alpha_j\in C^1(I_j)$, being $\dot\alpha_j$ constant on $I_j$ for every $j=1,\dots,N$. By
\begin{align*}
V(\phi)=&\lim_{k\to\infty}V_k(\phi)=\int_{M_k}\phi(x,T_xM_k)\,d\Ha^1
\\
=&\lim_{k\to \infty}\sum_{j=1}^N\int_{I_j}\phi\Big(\alpha_{j,k}(s),
Id-\frac{\dot\alpha_{j,k}(s)\otimes\dot\alpha_{j,k}(s)}{\vert\dot\alpha_{j,k}(s)\vert^2}\Big) \vert\dot\alpha_{j,k}(s)\vert\,ds
\\
=&\sum_{j=1}^N\int_{I_j}\phi\Big(\alpha_{j}(s),
Id-\frac{\dot\alpha_{j}(s)\otimes\dot\alpha_{j}(s)}{\vert\dot\alpha_{j}(s)\vert^2}\Big) \vert\dot\alpha_{j}(s)\vert\,ds
\end{align*}
we conclude that  \eqref{eq:cooper} holds.

In order to prove \eqref{eq:sheldon} we proceed by contradiction. Suppose, without loss of generality, that $T_1\neq T_2$, and $T_1\cap B_R,\,T_2\cap B_R\subset \mathrm{spt}(\mu_{V}\cap B_R)$,  and $T_1\cap T_2\cap B_R\neq\emptyset$. We can find $\alpha_{j_l,k}\in C^2(I_{j_l})$, parametrizing a connected components of $M_k\cap B_{3R/2}$, uniformly convergent to   a constant speed paramatrization $\alpha_{j_l}\in C^1(I_{J_l})$ of $T_l\cap B_R$ ($l=1,2$). Since $\vert\dot\alpha_{j_l,k}\vert$ is constant for $l=1,2$ and every $k\in\N$, by the uniform convergence of  $\alpha_{j_1,k},\,\alpha_{j_2,k}$ and by $T_1\cap T_2\cap B_R\neq\emptyset$ we can conclude that $\alpha_{j_1,k}(I_{j_1})\cap\alpha_{j_2,k}(I_{j_2})\neq\emptyset$ for every $k$ big enough. But this contradicts the embededdness assumption made  on $M_k$. Hence \eqref{eq:sheldon} holds too, and the proof is complete.
\endproof

\begin{lemma}\label{lem:strong-conv-eps}
Let $\tueps\in C^2(B_{2R})$ be such that
\begin{gather}\label{eq:strong-conv-eps-I}
0<\liminf_{\eps\to 0}\int_{B_{2R}}\frac{\eps}{2}\vert\nabla \tueps\vert^2+\frac{W(\tueps)}{\eps}\,dx\leq\limsup_{\eps\to 0}\int_{B_{2R}}\frac{\eps}{2}\vert\nabla \tueps\vert^2+\frac{W(\tueps)}{\eps}\,dx<+\infty, 
\\
 \lim_{\eps\to 0}\frac{1}{\eps}
\int_{B_{2R}}\Big\vert\eps\nabla^2\tueps-\frac{W^\prime(\tueps)}{\eps}\nu_{\tueps}\otimes\nu_{\tueps}\Big\vert^2dx=0.
\label{eq:strong-conv-eps-II}
\end{gather}
Then, being  $\Vepstueps$ as in \eqref{eq:ridi-culo}, up to a subsequence we have
$$
\lim_{\eps\to0}\Vepstueps=\widetilde V\text{ as varifolds in }\Omega,
$$
where $\widetilde V\in\mathbf{IV}_1(B_{2R})$
 is stationary and
  verifies \eqref{eq:cooper} and \eqref{eq:sheldon}.
\end{lemma}

\proof
We begin by selecting a subsequence (not relabeled) such that 
$$
0<\lim_{\eps\to 0} \int_{B_{2R}}\frac{\eps}{2}\vert\nabla \tueps\vert^2+\frac{W(\tueps)}{\eps}\,dx<+\infty.
$$
We fix $\Omega^\prime$ such that $B_R\subset\subset\Omega^\prime\subset\subset B_{2R}$.
 By 
Sard's Lemma and \cite[Lemma 7.1]{BeMu:09} we can find a  
subsequence $\{\Vepstuepsk\}_{k}$ and  a subset $J\subset 
[-1,1]$, with $\LL^1(J)=0$, such that for every $s\in [-1,1]\setminus J$,
\begin{gather*}
\{\tuepsk=s\}\text{ is a smooth embedded 
surface without boundary in }\Omega^\prime
\\
 \{\tuepsk=s\}\cap\{\nabla \tuepsk=0\}=\emptyset,
\\
\lim_{k\to\infty} \var(\{\tuepsk=s\}
,1)=\widetilde V~ \text{ as varifolds on }\Omega^\prime.
\end{gather*}

For every $x\in\Omega^\prime$ such that $\tuepsk(x)=s\in [-1,1]\setminus J$ we set
$$
\sff_{\tuepsk}:=\frac{(P^{\tuepsk})^T\nabla^2 \tuepsk P^{\tuepsk}}{\vert\nabla \tuepsk\vert}\otimes\nituepsk,
$$
that is $\sff_{\tuepsk}(x)$ is the second fundamental form of $\{\tuepsk=s\}$ at the point $x$ (see \eqref{eq:def-sff-lev-surf}). Let us also recall that (see \cite[Lemma 5.3]{BeMu:09})
\begin{equation}\label{equlo}
\vert\sff_{\tuepsk}\vert\, \eps_k\vert\nabla\tuepsk\vert\leq \Big\vert\eps\nabla^2\tueps-\frac{W^\prime(\tueps)}{\eps}\nu_{\tueps}\otimes\nu_{\tueps}\Big\vert.
\end{equation}

Next we fix $\delta>0$ and set $I_\delta:=[-1+\delta,1-\delta]$.  By \eqref{equlo}
we have
\begin{align*}
&\int_{I_\delta\setminus J}\left\vert\delta \var(\{\tuepsk=s\},1)\right\vert(\Omega^\prime)\,ds=
\int_{I_\delta\setminus J}\int_{\{\tuepsk=s\}\cap\Omega^\prime}
\left\vert \mathrm{div}
\left(
\nituepsk
\right)\right\vert\,d\Ha^1\,ds
\\
\leq & \frac{1}{(2\delta-\delta^2)}\int_{\Omega^\prime}\left\vert \mathrm{div}\left(
\nituepsk
\right)\right\vert\sqrt{2W(\tuepsk)}\vert\nabla \tuepsk\vert\,dx=
\frac{2}{(2\delta-\delta^2)}
\int_{\Omega^\prime}\vert\sff_{\tuepsk}\vert\sqrt{2W(\tuepsk)}\vert\nabla \tuepsk\vert\,dx
\\
\leq &\frac{2}{(2\delta-\delta^2)} \left(\int_{B_{2R}}\vert\sff_{\tuepsk}\vert^2\;d\widetilde\mu^\eps_{\tuepsk}\right)^{1/2}
\left(\int_{B_{2R}}\frac{W(\tuepsk)}{\eps_k}\,dx
\right)^{1/2}
\\
\leq &\frac{2}{(2\delta-\delta^2)} \left(\frac{1}{\eps_k}\int_{B_{2R}} \Big\vert\eps\nabla^2\tueps-\frac{W^\prime(\tueps)}{\eps}\nu_{\tueps}\otimes\nu_{\tueps}\Big\vert^2\;dx\right)^{1/2}
\left(\int_{B_{2R}}\frac{W(\tuepsk)}{\eps_k}\,dx
\right)^{1/2}.
\end{align*}
By the choice 
of  $\eps_k$ and of the set $J$,
and by
\eqref{eq:strong-conv-eps-II}, we can conclude that there exists 
$s_{\eps_k}\in I_\delta\setminus J$ such that, setting $V_k:=\var(\{\widetilde u_{\eps_k}=s_{\eps_k}\}\cap \Omega^\prime,1)$, we have
\begin{gather*}
\limsup_{k\to \infty}\mu_{V_k}(\Omega^\prime)<+\infty,\quad
\\
\limsup_{k\to \infty}\Big\vert\delta V_k\Big\vert(\Omega^\prime)=0.
\end{gather*}
therefore we are in a position to apply Lemma \ref{lem:app-smooth-flat} to the sequence $\{V_k\}_k\subset\mathbf{IV}_1(B_2R)$. Hence we can conclude the proof  by \cite[Theorem 4.1]{BeMu:09} and 
$$
\lim_{k\to\infty} \var(\{\widetilde u_{\eps_k}=s_{\eps_k}
\}
,1)=\widetilde V~ \text{ as varifolds on }\Omega^\prime.
$$
\endproof

\section{Proof of Theorem \ref{theo:approx}}\label{appotta}

Since $V$ is a Hutchinson's varifold with square-integrable second fundamental form by Remark \ref{rem:MC-Hutch-MC-Allard} the conclusions of Remark \ref{rem:briosch} hold.

For $x\in \R^2$ and $\lambda>0$ we define
\begin{gather*}
\eta_{x,\lambda}:\R^2\to\R^2,
\qquad
y\mapsto \frac{y-x}{\lambda}.
\end{gather*}
and  consider, for $x\in \mathrm{spt}(\mu_V)$, the Radon measure
$$
\mut_{x,\rho}(\psi):=\frac{1}{\rho}\int_{\eta_{x,\rho}(M)} \psi(y)\,\theta(\rho y+x)d\Ha^1(y), \quad \forall \psi\in C^0_c(\R^2).
$$
By \cite[Theorem 3.4]{Hutch-reg:84} (see also \cite{Hutch-reg:87}) we can conclude that for \textit{every} $x\in \mathrm{spt}(\mu_V)$  there exists a Radon measure $\mut_x$ on $\R^2$ such that
$$
\lim_{\rho\to 0^+}\mut_{x,\rho}(\psi)=\mut_x(\psi),\quad\forall \psi\in C^0_c(\Omega),
$$ 
and moreover that the measure $\mut_x$ satisfies
$$
\mut_x=\sum_{i=1}^{N_x}\Theta_{i}(x)\Ha^1_{\res \tilde T_i(x)},
$$
where $N_x\in\N$, $\tilde T_1(x),\dots,\tilde T_{N_x}(x)\in G_{1,2}$, and $\Theta_1(x),\dots,\Theta_{N_x}(x)\in\N$. In order to prove  the existence of an unique tangent line in every point of $\mathrm{spt}(\mu_V)$ we show that  $N_x= 1$ for every $x\in\mathrm{spt}(\mu_V)$. 

\noindent Without loss of generality we suppose that $x=0$. In view of the Hutchinson's regularity result cited above,  to conclude that $N_0=1$ it is enough to prove that for every sequence $\{\rho_k\}_k\subset\R^+$ such that $\lim_{k\to\infty}\rho_k=0$, setting 
$$
\mut_k(\psi):=\frac{1}{\rho_k}\int_{\eta_{0,\rho_k}(M)} \psi(y)\,\theta(\rho_ky)d\Ha^1(y), \quad \forall \psi\in C^0_c(\R^2),
$$
we have
\begin{gather}\label{aim}
\mut(\psi)=\Theta^1(\mu,0)\int_{T}\psi(y)\,d\Ha^1,
\end{gather}
where $T\in G_{1,2}$ is a linear $1$-dimensional subspace of
$\R^2$.

 Since $\mut_k\to\mut$ as Radon-measures on $\R^2$ and $\mu_{\eps_k}\to\mu_V$ as Radon measures in $\Omega$, for every open bounded  subset $U\subset\R^2$,  we can find a sequence $\{\eps_k\}_k$ such that 
\begin{gather*}
\lim_{k\to\infty}\eps_k=\lim_{k\to\infty}\frac{\eps_k}{\rho_k}=0,
\end{gather*}
and such that, setting $\tu_k(y):=u_{\eps_k}(\rho_ky)$, $\teps_k:=\eps_k/\rho_k$, the following hold
\begin{gather*}
\int_B\frac{\teps_k}{2}\vert\nabla \tu_k\vert^2+\frac{W(\tu_k)}{\teps_k}\,dx=\frac{\mu_{\eps_k}(\rho_kB)}{\rho_k}\to
\mut(B), ~\forall B\subset\subset U \text{ Borel},
\\
0<\lim_{k\to\infty}\int_U\teps_k\vert\nabla \tu_k\vert^2+\frac{W(\tu_k)}{\teps_k}\,dx<+\infty.
\end{gather*}
Moreover by the definition of $\tu_k$ and $\teps_k$ and \eqref{eq:astianatte}, we have
\begin{gather}
\frac{1}{\teps_k}\int_U \Big \vert \teps_k\nabla^2 \tu_k
-\frac{W^\prime(\tu_k)}{\teps_k}
\nu_{\tu_k}
\otimes\nu_{\tu_k}
\Big\vert^2\,dy
\notag
\\
=\frac{\rho_k}{\eps_k}
\int_U\left\vert\eps_k\nabla^2 u_{\eps_k}(\rho_k y)-\frac{W^\prime(u_{\eps_k}(\rho_k y))}{\eps_k}\frac{\nabla u_{\eps_k}(\rho_k y)\otimes
\nabla u_{\eps_k}(\rho_k y)}{\vert\nabla u_{\eps_k}(\rho_k y)\vert^2}\right\vert^2\rho_k^2\,dy
\\
=\frac{\rho_k}{\eps_k}\int_{\rho_kU}\Big\vert
\eps_k\nabla^2 u_{\eps_k}(x)-\frac{W^\prime(u_{\eps_k}(x))}{\eps_k}
\nu_{u_{\eps_k}}(x)\otimes\nu_{u_{\eps_k}}(x)\Big\vert^2\,dx\leq C\rho_k.
\notag
\end{gather}
We can thus apply Lemma \ref{lem:strong-conv-eps} and extract a sequence (not relabelled) such that 
\begin{gather*}
V^{\teps_k}_{\tu_k}\to\widetilde V=\sum_{j=1}^N\Theta_j\Ha^1_{\res T_j\cap U}
 \text{ as varifolds in }U,
\\
\mu_{\widetilde V}=\mut_{\res U},
\end{gather*}
and $T_i\cap T_j\cap B_R=\emptyset$ for every $B_{2R}\subset\subset U$.

However since $\mut$ verifies
$$
\mut(B_R)=R\lim_{k\to\infty}\frac{\mut_k(B_{\rho_kR})}{R\rho_k}=R\,\Theta(\mu_V,0)=:R\,\Theta_0,
$$
we have $\widetilde V=\Theta_0\Ha^1_{\res T\cap U}$ that is \eqref{aim}. 
\qed

\section{Proof of Corollary \ref{cor:sheldon-brown}}

We begin proving the so-called $\Gamma-\liminf$-inequality.
We suppose that $\{u_\eps\}_\eps\subset X$ satisfies \eqref{eq:mio} (otherwise we have nothing to prove). By Theorem \ref{the:mio} we can find a subsequence $\{\eps_k\}_{k\in\N}$  such that $\lim_{k\to\infty}\eps_k=0$ and
\begin{gather*}
\lim_{k\to \infty}\mathcal E_{\eps_k}(u_{\eps_k})=\liminf_{\eps\to 0}\mathcal E_\eps(u_\eps),
\\
L^1(\Omega)-\lim_{\eps\to 0}u_{\eps_k}\to u=\ind_E \in BV(\Omega,\{-1,1\}),
\\
 \lim_{k\to\infty}V^{\eps_k}_{u_{\eps_k}}=V\in\mathscr{CV}^2_1(\Omega)\text{ as varifolds}.
\end{gather*}
If we prove that $V\in\mathscr D$, by Theorem \ref{theo:rel-el} and Theorem \ref{the:mio}-(C1), we obtain
\begin{gather*}
\liminf_{\eps\to 0}\mathcal E_\eps(u_\eps)
=\lim_{k\to\infty}\mathcal E_{\eps_k}(u_{\eps_k})
\geq
\surftens\int(1+\vert \Hv_V\vert^2)\,d\mu_V=\mathcal F_2(V)\geq \surftens\overline{\mathcal F}(E).
\end{gather*}
That is the $\Gamma-\liminf$ inequality holds. In order to prove that $V=\var(M,\theta)\in\mathscr{D}(\Omega)$, by  
\eqref{eq:Cahr-D},  it is enough  to show that: (i) $\var(M,\theta)$ is actually a Hutchinson's curvature varifold in the whole of $\R^2$; (ii) $M\cup\partial\Omega$ has an unique tangent-line in \textit{every} point.

We begin establishing that  (i) holds. To this purpose we fix  $\Omega_1\subset\subset\R^2$ such that $\Omega_1\supset\supset\Omega$, and define $\Omega_\delta:=\{x\in\Omega_1\setminus\Omega:~\dist(x,\Omega)<\delta\}$ for $\delta>0$. Next
we  notice that, since  $\{u_\eps\}_\eps\subset X$,
 we can extend  $u_\eps$ to $u^\prime_\eps\in W^{2,2}(\Omega_1)$ simply setting $u^\prime_\eps\equiv 1$ on $\Omega_1\setminus\overline\Omega$. 
Since $u_\eps^\prime$ satisfies \eqref{eq:mio} on $\Omega_1$, by Theorem \ref{the:mio} we can extract a subsequence such that $\lim_{\eps\to 0}V^\eps_{u^\prime_\eps}=V^\prime\in\mathscr{CV}^2_1(\Omega_1)$. However, since 
$$
\int_{\Omega^\prime}\frac{\eps}{2}\vert\nabla u^\prime_\eps\vert^2+\frac{W(u^\prime_\eps)}{\eps}\,dx=\int_{\Omega}\frac{\eps}{2}\vert\nabla u_\eps\vert^2+\frac{W(u_\eps)}{\eps}\,dx,
$$
we can conclude that $\mathrm{spt}(\mu_{V^\prime})\subset\overline\Omega$. Hence we obtain that $V=V^\prime$ is a Hutchinson varifold in $\Omega_1$ whose support is compactly contained in $\Omega_1$, and therefore $\var(M,\theta)\in\mathscr{CV}^2_1(\R^2)$.

We now pass to prove (ii). By \cite[Theorem 4.2]{BeMu:09}, we can find an infinitesimal, strictly decreasing sequence $\{\delta_\eps\}_\eps\subset\R^+$, and a sequence $g_\eps\in C^2(\Omega_1\setminus\overline\Omega)$ such that 
\begin{gather*}
g_\eps\equiv -1 \text{ on }\Omega_1\setminus\Omega_{\delta_\eps}, \quad g_\eps\equiv 1 \text{on }\Omega_{2\eps}\setminus\Omega,
\\
\lim_{\eps\to 0}V^\eps_{g_\eps}=\var(\partial\Omega,1)\text{ as varifolds },\quad L^1(\Omega_1\setminus\overline\Omega)-\lim_{\eps\to 0}g_\eps\equiv -1,
\\
\int_{\Omega_1\setminus\overline\Omega}\frac{\eps}{2}\vert\nabla g_\eps\vert^2+\frac{W(g_\eps)}{\eps}\,dx
=\surftens\Ha^1(\partial \Omega)+O(\eps),
\\
\frac{1}{\eps}\int_{\Omega_1\setminus\overline\Omega}\left\vert\eps\nabla^2 g_\eps
-\frac{W^\prime(g_\eps)}{\eps}\nu_{g_\eps}\otimes\nu_{g_\eps}\right\vert^2\,dx =\int_{\partial\Omega}\vert \sff_{\partial\Omega}\vert^2\,d\Ha^1+O(\eps).
\end{gather*}

Hence, again by the assumption $\{u_\eps\}\subset X$, we can conclude that setting
\begin{gather*}
u^{\prime\prime}_\eps(x):=\begin{cases}
                                u_\eps(x)  &\text{if }x\in\Omega,
                                \\
                                g_\eps               &\text{if }x\in \Omega_{1}\setminus\overline\Omega,
                               \end{cases}
\end{gather*}
the sequence $\{u_\eps^{\prime\prime}\}_\eps\subset W^{2,2}(\Omega_1)$ satisfies \eqref{eq:mio}, and moreover (up to subsequences)  as $\eps\to 0$, we have
$$
V^\eps_{u^{\prime\prime}_\eps}\to \var(M,\theta)+\var(\partial\Omega,1)\in\mathscr{CV}^2_1(\R^2)\text{ as varifolds}.
$$  
Applying Theorem \ref{theo:approx} to the sequence $\{u^{\prime\prime}_\eps\}_\eps$ we obtain that $M\cup\partial \Omega$ has an unique tangent line in every point. Hence $V\in\mathscr D(\Omega)$ and the $\Gamma-\liminf$ inequality holds.

Finally, the $\Gamma-\limsup$ inequality now follows by \cite[Theorem 4.2]{BeMu:09} and a standard density argument. In fact, by the previous step we can conclude that for every $u=\ind_E\in L^1(\Omega)$ such that $\Galim\mathcal E_{\eps,\res X}(u)<+\infty$ we also have $\overline{\mathcal F}(u)<+\infty$, and therefore we can find a sequence $\{E_h\}_h$ such that
$E_h\subset\subset\Omega$, and $\Omega\cap\partial E_h\in C^2$, and
$$
L^1(\Omega)-\lim_{h\to\infty}\ind_{E_h}=u,\quad \lim_{h\to\infty}\mathcal F(\ind_{E_h})=\overline{\mathcal F}(u).
$$
\qed

\section{Proof of Theorem \ref{theo:sega}}\label{smoothie}

Without loss of generality we suppose that $\Omega=B_1$.
In order to prove Theorem \ref{theo:sega} we begin by
 showing
 the existence of a sequence
 $\{u_\eps\}_\eps\subset C^3(\Omega)$ such that 
\begin{gather}
\eps\Delta u_\eps-\frac{W^\prime(u_\eps)}{\eps}=0,\qquad\forall\eps>0
\label{eq:scajo}
\\
\sup_{\eps>0}\int_\Omega\frac{\eps}{2}\vert\nabla u_\eps\vert^2+\frac{W(u_\eps)}{\eps} \,dx\leq C,
\label{eq:ladro}
\end{gather}
and 
\begin{itemize}
\item[(a)] $L^1(\Omega)-\lim_{\eps\to 0}u_\eps=u=\ind_{E_0}\in BV(\Omega,\{-1,1\})$, where
$$
E_0=\{(x_1,x_2)\in\Omega:~ x_1> 0,\,x_2> 0\}\cup \{(x_1,x_2)\in\Omega:~ x_1< 0,\,x_2 <0\};
$$
\item[(b)] $\lim_{\eps\to 0}V_\eps=V_{\mathfrak C}:=\var(\mathfrak C\cap \Omega,\theta)\in\mathbf{IV}_1(\Omega)$ as varifolds, where 
$$
\mathfrak C:=\{x:=(x_1,x_2)\in\R^2:~x_1\text{ vel }x_2 \text{ equals }0\}=\partial E_0.
$$  
\end{itemize}
 The fact
that showing the existence of a sequence $\{u_\eps\}_\eps$ with the above properties is enough to conclude the proof of the first part of
Theorem \ref{theo:sega} is pretty easy to see. In fact, 
 since $0\in \mathfrak C=
\mathrm{spt}(\mu_{V_{\mathfrak C}})$ and the tangent cone 
in $0$ to $\mathfrak C$ coincides with
 $\mathfrak C$ itself
,  we have that $\mu_{V_{\mathfrak C}}$ can not have an uniquely defined tangent line in  $0\in\mathrm{spt}(\mu_{V_{\mathfrak C}})\cap\Omega$.

We construct the sequence $\{u_\eps\}\subset C^3(\Omega)$ verifying \eqref{eq:scajo}, \eqref{eq:ladro} via the blow-down of a particular entire solution of the Allen-Cahn equation in the plane. More precisely,
let $U\in C^3(\R^2)$ be a ``saddle solution'' of the Allen-Cahn equation, that is
\begin{equation}\label{eq:AC}
\Delta U=W^\prime(U)\qquad \text{on }\R^2,
\end{equation}
and $U$ is such that
\begin{itemize}
\item $\|U\|_{L^\infty(\R^2)}\leq 1$,  $\{U=0\}=\mathfrak C$ and $U>0$ (respectively $U<0$) in the I and III (respectively II and IV) quadrant of $\R^2$;
\item there exists $C>0$ such that for every $R>0$ 
\begin{gather}\label{eq:Xavi}
\int_{B_R}\frac{1}{2}\vert\nabla U\vert^2+W(U)\,dy\leq C\,R.
\end{gather}
\end{itemize}
The existence of such a solution has been proved in  \cite[Theorem 1.3]{CabTe:JEMS} (see also \cite{Dang-Fife-Pel, DelPino}).

We define $\{u_\eps\}_\eps\subset C^2(\Omega)$  by $u_\eps(x):=U(x/\eps)$. By \eqref{eq:AC}, \eqref{eq:Xavi} we then have
\begin{gather*}
\eps\Delta u_\eps-\frac{W^\prime(u_\eps)}{\eps}=0,
\\
\int_\Omega\frac{\eps}{2}\vert\nabla u_\eps\vert^2+\frac{W(u_\eps)}{\eps} \,dx=\eps\int_{B_{\eps^{-1}}}\frac{1}{2}\vert\nabla U\vert^2+W(U)\,dy\leq C,
\end{gather*}
that is \eqref{eq:scajo} and \eqref{eq:ladro} hold.
Hence we are in a position to apply the results proved in \cite{HuTone}, and obtain that
\begin{itemize}
\item[(HT1)] (see \cite[Proposition 2.2]{HuTone}) for every $r<1$ there exists $c:=c(r)>0$ such that $\sup_{B_r}\xi^+_\eps\leq c$ for every $\eps$ small enough ;
\item[(HT2)] (see \cite[Proposition 3.4]{HuTone}) for every $x\in\Omega$, $0<\sigma<\rho$ such that $B_\rho(x)\subset\subset B_r$ ($r<1$), and $\eps$ small enough we have
\begin{gather}\label{eq:monoton}
\frac{\mu_\eps(B_\rho(x))}{\rho}\geq \frac{\mu_\eps(B_\sigma(x))}{\sigma}-c\rho,
\end{gather}
where $c=c(r)$ is defined in (HT1);
\item[(HT3)] (see \cite[Theorem 1]{HuTone}) from the sequence $\{V_\eps\}_\eps$ (see \eqref{eq:ridi-culo})  we can extract a subsequence (not relabeled) such that 
$$
\lim_{\eps\to 0}V_\eps= V:=\var(M,\theta) \text{ as varifolds in }\Omega,
$$ 
and $V\in\mathbf{IV}_1(\Omega)$ is stationary.
\end{itemize}
Next we show that $\mu_V(\Omega)>0$ and $M=\mathrm{spt}(\mu_V)= \mathfrak C\cap\Omega$.

Let $x_0\in\mathfrak C\cap \Omega$. We choose $\eps$ small enough that $B_\eps(x_0)\subset\subset B_{(1-\vert x_0\vert)/2}\subset\subset\Omega$, and define
\begin{gather*}
\tilde U_\eps\in C^2(B_1),\quad 
\tilde U_\eps(z):=u_\eps(\eps z+x_0)=U(z+\eps^{-1}x_0).
\end{gather*}
We then have
\begin{gather*}
\Delta\tilde U_\eps=W^\prime(\tilde U_\eps)\text{ in }B_1,
\text{ and } \tilde U_\eps(0)=U(\eps^{-1}x_0)=0.
\end{gather*}
Hence, by standard elliptic estimates, we have $\|\tilde U_\eps\|_{C^1(B_{1/2})}<\tilde C$, where $\tilde C>0$  is uniform with respect to $\eps$, and therefore we can find $\delta>0$ (independent of $\eps$) such that 
$\sup_{z\in B_\delta}\vert\tilde U_\eps(z)\vert<1/2$. Hence
\begin{align*}
\frac{\mu_\eps(B_{\delta\eps}(x_0))}{\delta\eps}=&\frac{1}{\delta\eps}\int_{B_{\delta\eps}(x_0)}\frac{\eps}{2}
\vert\nabla u_\eps\vert^2+\frac{W(u_\eps)}{\eps}\,dx
\\
=&
\int_{B_\delta}\frac{1}{2}\vert\nabla \tilde U_\eps\vert^2+W(\tilde U_\eps)\,dz
\geq \int_{B_\delta}W(\tilde U_\eps)\,dz\geq C_W
\end{align*}
where $C_W:=C_W(\delta)=\pi\delta^2\min\{W(s):\,s\in (-1/2,1/2)\}$.

We now choose $\rho<\rho_0$ where $\rho_0$ is such that $C_W-c\rho_0>C_W/2$. By \eqref{eq:monoton}  for every $\eps$ small enough we have
$$
\frac{\mu_\eps(B_\rho(x_0))}{\rho}\geq \frac{\mu_\eps(B_{\delta\eps}(x_0))}{\delta\eps}-c\rho>\frac{C_W}{2},
$$
from which we deduce $\mu_V(\Omega)>0$ and $\mathrm{spt}(\mu_V)\supseteq\mathfrak C\cap \Omega$.

However, in view of  \cite[Lemma 5]{Dang-Fife-Pel}, we can find a  constant $K>0$, independent of $\eps$, such that
that for every $\eta\in(0,1/2)$ there exists $\eps_0:=\eps_0(\eta)$ such that for $\eps<\eps_0$ we have
$$
\{x=(x_1,x_2)\in\Omega:~\vert x_1\vert,\,\vert x_2\vert>\eps k\}\subset\{x\in\Omega:~\vert u_\eps(x)\vert \geq 1-\eta\}.
$$
By this latter estimate and \cite[Proposition 5.1]{HuTone}, we can conclude that 
$$
\lim_{\eps\to 0}\mu_\eps(\overline A)=0,\quad \forall A \subset\subset \Omega\setminus\mathfrak C.
$$
Hence $\mathrm{spt}(\mu_V)\subseteq\mathfrak C$ and this concludes the proof of the part of Theorem \ref{theo:sega}.

It remains to prove that \eqref{doggy} holds. To this aim it is enough to remark that, being $\{u_\eps\}_\eps$ and $u$ as above,
by Proposition \ref{lem:labradford} and Theorem \ref{theo:approx} we have
$$
\overline{\mathcal F}_o(u)=\Galim \Eeps(u)=+\infty.
$$

\qed
\bibliography{rel}

\def\cprime{$'$} \def\cprime{$'$}
\begin{thebibliography}{10}

\bibitem{BDMP}
G.~Bellettini, G.~Dal~Maso, and M.~Paolini.
\newblock Semicontinuity and relaxation properties of a curvature depending
  functional in $2$d.
\newblock {\em Ann. Scuola Norm. Sup. Pisa Cl. Sci.}, 20(2):247--297, 1993.

\bibitem{BeMu:07}
G.~Bellettini and L.~Mugnai.
\newblock A varifolds representation of the relaxed elastica functional.
\newblock {\em J. Convex Anal.}, 14(3):543--564, 2007.

\bibitem{BeMu:09}
G.~Bellettini and L.~Mugnai.
\newblock Approximation of the {H}elfrich's functional via diffuse interfaces.
\newblock to appear on \textit{SIAM J. Math. Anal}, 2010.

\bibitem{BePa:93}
G.~Bellettini and M.~Paolini.
\newblock Approssimazione variazionale di funzionali con curvatura.
\newblock {\em Seminario Analisi Matematica Univ. Bologna}, 1993.

\bibitem{BraMar}
A.~Braides and R.~March.
\newblock Approximation by {$\Gamma$}-convergence of a curvature-depending
  functional in visual reconstruction.
\newblock {\em Comm. Pure Appl. Math.}, 59(1):71--121, 2006.

\bibitem{CabTe:JEMS}
X.~Cabr\'e and J.~Terra.
\newblock Saddle-shaped solutions of bistable diffusion equations in all of
  $\mathbb{R}^{2m}$.
\newblock {\em JEMS}, 43:819--943, 2009.

\bibitem{CampHern2}
F.~Campelo and A.~Hernandez-Machado.
\newblock Shape instabilities in vesicles: A phase-field model.
\newblock {\em The European Physical Journal}, 143(1):101--108, 2007.

\bibitem{Suka}
T.~Chan, S.~Kang, and J.~Shen.
\newblock Euler's elastica and curvature-based inpainting.
\newblock {\em SIAM J. Appl. Math.}, 63(2):564--592, 2002.

\bibitem{DM}
G.~Dal~Maso.
\newblock {\em An introduction to {$\Gamma$}-convergence}, volume~8 of {\em
  Progress in Nonlinear Differential Equations and their Applications}.
\newblock Birkh\"auser, Boston, MA, 1993.

\bibitem{Dang-Fife-Pel}
H.~Dang, P.~Fife, and L.~Peletier.
\newblock Saddle solutions of the bistable diffusion equation.
\newblock {\em Z. Angew. Math. Phys.}, 43:984--998, 1992.

\bibitem{DG}
E.~De~Giorgi.
\newblock Some remarks on {$\Gamma$}-convergence and least squares method.
\newblock In {\em Composite media and homogenization theory (Trieste, 1990)},
  volume~5 of {\em Progr. Nonlinear Differential Equations Appl.}, pages
  135--142. Birkh\"auser Boston, Boston, MA, 1991.

\bibitem{DelPino}
M.~del Pino, M.~Kowalczyk, F.~Pacard, and J.~Wei.
\newblock Multiple-end solutions to the {A}llen-{C}ahn equation in $\mathbb{
  R}^2$.
\newblock {\em J. Funct. Anal}, 258(2):458--503, 2010.

\bibitem{MeMaPa}
P.~Dondl, L.~Mugnai, and M.~R\"oger.
\newblock Confined elastic curves.
\newblock preprint, 2010.

\bibitem{DuWill}
Q.~Du, C.~Liu, R.~Ryham, and X.~Wang.
\newblock A phase field formulation of the {W}illmore problem.
\newblock {\em Nonlinearity}, 18(3):1249--1267, 2005.

\bibitem{DuCaz}
Q.~Du, C.~Liu, R.~Ryham, and X.~Wang.
\newblock Diffuse interface energies capturing the {E}uler number: relaxation
  and renormalization.
\newblock {\em Commun. Math. Sci.}, 8(1):233--242, 2007.

\bibitem{DuUno}
Q.~Du, C.~Liu, and X.~Wang.
\newblock A phase field approach in the numerical study of the elastic bending
  energy for vesicle membranes.
\newblock {\em J. Comput. Phys.}, 198(2):450--468, 2004.

\bibitem{EseShen}
S.~Esedoglu and J.~Shen.
\newblock Digital inpainting based on the mumford-shah-euler image model.
\newblock {\em European J. Appl. Math.}, 13(4):353--370, 2002.

\bibitem{Hutch-reg:84}
J.~Hutchinson.
\newblock ${C}^{1,\alpha}$-multiple function regularity and tangent cone
  behaviour for varifolds with second fundamental form in ${L}^p$.
\newblock In {\em Geometric measure theory and the calculus of variations
  (Arcata, Calif., 1984)}, volume~44 of {\em Proc. Sympos. Pure Math.}, pages
  281--306. Amer. Math. Soc., Providence, RI, 1984.

\bibitem{Hu}
J.~Hutchinson.
\newblock Second fundamental form for varifolds and the existence of surfaces
  minimising curvature.
\newblock {\em Indiana Univ. Math. J.}, 35(1):281--306, 1986.

\bibitem{Hutch-reg:87}
J.~Hutchinson.
\newblock Some regularity theory for curvature varifolds. miniconference on
  geometry and partial differential equations.
\newblock In {\em {M}iniconference on geometry and partial differential
  equations ({C}anberra, {J}une, 1986)}, volume~12 of {\em Proc. Centre Math.
  Anal. Austral. Nat. Univ}, pages 60--66. Austral. Nat. Univ, Canberra, 1987.

\bibitem{HuTone}
J.~Hutchinson and Y.~Tonegawa.
\newblock Convergence of phase interfaces in the van der
  {W}aals-{C}ahn-{H}illiard theory.
\newblock {\em Calc. Var. Partial Differential Equations}, 10(1):49--84, 2000.

\bibitem{KORV}
R.~Kohn, F.~Otto, M.~Reznikoff, and E.~Vanden-Eijnden.
\newblock Action minimization and sharp-interface limits for the stochastic
  {A}llen-{C}ahn equation.
\newblock {\em Comm. Pure Appl. Math.}, 60(3):393--438, 2007.

\bibitem{Masnou}
G.~Leonardi and S.~Masnou.
\newblock Locality of the mean curvature of rectifiable varifolds.
\newblock {\em Adv. Calc. Var.}, 2(1):17--42, 2009.

\bibitem{LiuSatoTone}
C.~Liu, N.~Sato, and Y.~Tonegawa.
\newblock On the existence of mean curvature flow with transport term.
\newblock to appear on Interfaces Free Bound., 2009.

\bibitem{LoMar}
P.~Loreti and R.~March.
\newblock Propagation of fronts in a nonlinear fourth order equation.
\newblock {\em European J. Appl. Math.}, 11(2):203--213, 2000.

\bibitem{Lowe}
J.~S. Lowengrub, A.~R\"atz, and A.~Voigt.
\newblock Phase-field modeling of the dynamics of multicomponent vesicles:
  spinodal decomposition, coarsening, budding, and fission.
\newblock {\em Phys. Rev. E}, 79(3):82C99--92C10, 2009.

\bibitem{Men}
U.~Menne.
\newblock Second order rectifiability of integral varifolds of locally bounded
  first variation.
\newblock preprint, 2008.

\bibitem{MM}
L.~Modica and S.~Mortola.
\newblock Un esempio di {$\Gamma ^{-}$}-convergenza.
\newblock {\em Boll. Un. Mat. Ital. B}, 14(5):285--299, 1977.

\bibitem{MuRo2}
L.~Mugnai and M.~R\"oger.
\newblock Convergence of perturbed allen-cahn equations to forced mean
  curvature flow.
\newblock to appear on \textit{Indiana Univ. Math. J.}

\bibitem{MuRo1}
L.~Mugnai and M.~R\"oger.
\newblock The {A}llen-{C}ahn action functional in higher dimensions.
\newblock {\em Interfaces Free Bound.}, 10(1):45--78, 2008.

\bibitem{RoSch}
M.~R\"oger and R.~Sch\"atzle.
\newblock On a modified conjecture of {D}e {G}iorgi.
\newblock {\em Mathematische Zeitschrift}, 254(4):675--714, 2006.

\bibitem{Sato}
N.~Sato.
\newblock A simple proof of convergence of the {A}llen-{C}ahn equation to
  {B}rakke's motion by mean curvature.
\newblock {\em Indiana Univ. Math. J.}, 57(4):1743--1751, 2008.

\bibitem{Sch:09}
R.~Sch\"atzle.
\newblock Lower semicontinuity of the {W}illmore functional for currents.
\newblock {\em J. Differential Geom.}, 81(2):437--456, 2009.

\bibitem{Serfa}
S.~Serfaty.
\newblock {G}amma-convergence of gradient flows on hilbert and metric spaces
  and applications.
\newblock preprint, 2010.

\bibitem{Sim}
L.~Simon.
\newblock {\em Lectures on Geometric Measure Theory}, volume~3 of {\em
  Proceedings of the Centre for Mathematical Analysis, Australian National
  University}.
\newblock Centre for Math. Anal. Australian National Univ., Canberra, 1984.

\bibitem{ToneYuko}
Y.~Tonegawa and Y.~Nagase.
\newblock A singular perturbation problem with integral curvature bound.
\newblock {\em Hiroshima Math. Journal}, 37(3):455--489, 2007.

\end{thebibliography}
\bibliographystyle{abbrv}

\end{document}